\documentclass{svjour3}       

\smartqed  

\usepackage{moreverb}

\newcommand{\subfigimg}[3][,]{%
  \setbox1=\hbox{\includegraphics[#1]{#3}}
  \leavevmode\rlap{\usebox1}
  \rlap{\hspace*{25pt}\raisebox{\dimexpr\ht1-0.8\baselineskip}{#2}}
  \phantom{\usebox2}
}

\usepackage{amsmath, bm}
\usepackage{caption}
\usepackage{booktabs}

\usepackage{hyperref}
\usepackage{color}
\usepackage{multirow}

\usepackage{amssymb}
\usepackage{latexsym}
\usepackage{fancyhdr}
\usepackage{graphicx}
\usepackage{newlfont}
\usepackage{textcomp}
\usepackage{float}
\usepackage{multirow}

\usepackage{enumitem}

\newtheorem{assumption}{Assumption}

\usepackage{color}
\newcommand{\Giacomo}[1]{\textcolor{black}{{#1}}}

\begin{document}

\title{Approximation of probability density functions for {PDEs with random parameters} using truncated series expansions
\thanks{MG and HW thank the Isaac Newton Institute for Mathematical Sciences, Cambridge, for support and hospitality during the program {\em Uncertainty Quantification for Complex Systems: Theory and Methodologies} where work on this paper was undertaken. This work was supported by UK EPSRC grant numbers EP/K032208/1 and EP/R014604/1. GC and MG were supported in part by the US Air Force Office of Scientific Research grant FA9550-15-1-0001 and by the US Department of Energy Office of Science grant DE-SC0016591.}
}

\titlerunning{Estimation of PDFs for PDEs with random parameters using truncated series}        

\author{Giacomo Capodaglio \and Max Gunzburger \and Henry P. Wynn}

\authorrunning{G. Capodaglio, M. Gunzburger, and H. Wynn}

\institute{Giacomo Capodaglio \at
Department of Scientific Computing, Florida State University,  Tallahassee, FL 32306, USA.            
{Current address:} Computational Physics and Methods Group - Los Alamos National Laboratory, Los Alamos, New Mexico, USA; \email{gcapodaglio@lanl.gov}
           \and
Max Gunzburger \at
Department of Scientific Computing, Florida State University,  Tallahassee, FL 32306, USA; \email{mgunzburger@fsu.edu}
              \and
Henry P. Wynn \at 
The London School of Economics and Political Science, London WC2A 2AE, UK; \email{h.wynn@lse.ac.uk}
}

\dedication{Dedicated to Professor Enrique Zuazua on the occasion of his 60th birthday}

\date{Received: March 20, 2020 / Accepted: September 23, 2020}

\maketitle

\begin{abstract}
The probability density function (PDF) of a random variable associated with the solution of a partial differential equation {(PDE) with random parameters} is approximated using a truncated series expansion.
The {random PDE} is solved using two stochastic finite element methods, Monte Carlo sampling and the stochastic Galerkin method with global polynomials. 
The random variable is a functional of the solution of the {random PDE}, such as the average over the physical domain. The truncated series are obtained considering a finite number of terms in the Gram-Charlier or Edgeworth series expansions.
These expansions approximate the PDF of a random variable in terms of another PDF, and involve coefficients that are functions of the known cumulants of the random variable.
To the best of our knowledge, their use in the framework of {PDEs with random parameters} has not yet been explored.
\end{abstract}


\section{Introduction}
In the field of uncertainty quantification, the problem of describing the probability density function (PDF) of the output of a model given a known distribution of the input is of major relevance \cite{chen2016model}.
For instance, for the solution of inverse problems involving partial differential equations (PDEs) with Bayesian inference \cite{stuart2010inverse,dashti2017bayesian,bui2013computational}, the PDF for the forward problem may be sought to get insights on the forward uncertainty propagation \cite{marzouk2007stochastic,ma2009efficient}.
{In such a framework, the forward PDF is usually inferred from the histogram}, which may be obtained according to two different strategies.
The first is a direct approach, where Monte Carlo (MC) sampling is used to obtain realizations of the quantity of interest (QoI) by solving the forward problem for each sample.
Alternatively, once an approximation of the stochastic solution of the forward problem is obtained,
realizations may be obtained with MC sampling without solving the problem for each sample.
This is the method used for instance in \cite{marzouk2007stochastic}, where the stochastic solution of the forward problem is approximated with polynomial chaos, or in \cite{ma2009efficient}, where adaptive sparse grid collocation is employed.
{It is important to observe that despite providing a simple and easy approximation {of a density function}, histograms are not continuous, hence they cannot be differentiated if necessary. The differentiability of the estimator is an important feature in many applications \cite{cheng1995mean,fukunaga1975estimation,chacon2013data,sasaki2015direct,noh2014bias}.}
{Another popular density estimation technique that provides a smooth alternative to the histogram is the method of kernel density estimation (KDE) \cite{izenman1991review}. The performance of the KDE is dependent on the choice of a parameter called the bandwidth, which regulates the smoothness of the estimator. The kernels, which obviously also influence the smoothness, are usually centered at the samples, but a binning of the sample set can be performed although usually introducing an extra dependence on the binning size.}

In this work, we propose an alternative strategy to the histogram and the KDE for the estimation of the PDF of a quantity of interest, given a known distribution for the input data. The quantities of interest are associated with the solution of a {PDE with random parameters}, which is approximated using a stochastic finite element {(SFE)} method.
{The SFE methods used in this work are Monte Carlo sampling and the stochastic Galerkin (SG) method with global polynomials.}
Once the QoI has been obtained from the {solution of the random PDE using a SFE method}, its PDF is estimated by means of a truncated Gram-Charlier (GC) or Edgeworth (ED) series expansion \cite{kendall1943advanced}. {The use of GC and ED truncated series is the alternative approach we are considering.} These series have been used in several fields such as chemistry \cite{olive1991gram,di2001mathematical}, finance \cite{jondeau2001gram,popovic2012easy,niguez2012forecasting}, physics and astrophysics \cite{eggers2011determining,de2011edgeworth,o1992using,juszkiewicz1993weakly,blinnikov1998expansions,contaldi1999photographing}, material science \cite{rickman2015calculating}, oceanology \cite{zapevalov2011simulating}, power systems engineering \cite{fan2012probabilistic}, and other branches of applied mathematics \cite{pender2014Gram}.
However, to the best of our knowledge, their use in the field of {random PDEs} has not yet been explored.
The GC and ED expansions approximate a PDF adding successive corrections to a known PDF that is used as a first approximation.
The known PDF is usually chosen to be Gaussian, and only a few efforts have been made to generalize the GC expansion to the case of a non-Gaussian kernel \cite{brenn2017revisit,berberan2007expressing}. 
The terms in the GC and ED series involve derivatives of the known input PDF and are functions of the cumulants of the output PDF to be approximated.
The cumulants of a given random variable can be obtained analytically from its moments, which involve integrals of powers of the random variable \cite{kendall1943advanced}. {The advantage of employing the SG method with global polynomials for the solution of the random PDE is that with this choice, QoIs can be defined as polynomial functions in the random variable, so exact moments can be computed with appropriate quadrature rules. Consequently, also the coefficients of the GC and ED expansions can be computed exactly.}
The {SG} method is validated against the Monte Carlo method for small input variances of the random data.
\Giacomo{The approach we propose still makes use of the histogram to determine the PDF of the QoI, although the histogram is used only as a crude approximation to determine the most appropriate truncation order of the series. The approximate PDF is given analytically as a truncated series and it does not depend directly on the histogram, so it is not affected by the dimension of the bins used to construct it and by the number of samples. {Moreover, it does not depend on a smoothing parameter, such as the bandwidth, as does the KDE}. 
Finally, if the kernel is a $C^{\infty}$ function (i.e. a Gaussian distribution), the truncated GC and ED expansions are continuous and infinitely differentiable, therefore potentially more useful than just a histogram approximation {and simpler to deal with, given that smoothness is achieved without having to deal with extra parameters as in the KDE}.}

The paper is structured as follows: in Section \ref{formulation}, the mathematical problem and the input data description are introduced; in Section \ref{SFEMs}, the two {SFE methods}  employed for the solution of a {random PDE} are described, namely the Monte Carlo method and the Stochastic Galerkin method; the theory on the Gram-Charlier and Edgeworth expansions is laid out in Section \ref{expansions}. Arguments on asymptotic expansions and convergence are given in Section \ref{asymptexp} and Section \ref{convconsiderations}, respectively. Numerical results are reported in Section \ref{numResults}, for different types of output distributions.
The paper is concluded with Section \ref{discussion}, where our findings are discussed.

\section{Formulation of the problem}\label{formulation}
For {PDEs with random parameters}, the stochasticity is taken into account assuming that the input data depends on a random variable, other than the physical variable as in standard partial differential equations \cite{tartakovsky2011pdf,cliffe2011multilevel,nobile2008sparse,xiu2005high,nobile2008anisotropic,wan2005adaptive}.
The input data consists of coefficient functions and forcing term. For simplicity, here it is assumed that the stochastic contribution is introduced only by the coefficients and not by the forcing term.
The {random model} considered is the Poisson problem defined on $D \times \Omega$, where $D \subset \mathbb{R}^d$ and $\Omega$ is the sample space
\begin{align}\label{poissonSys}
\begin{cases}
- \nabla \cdot (\,a(\mathbf{x}, \omega) \,\, \nabla u(\mathbf{x}, \omega) \,)  = f(\mathbf{x}) & \text{in} \quad D \times \Omega \\
 u(\mathbf{x}, \omega)=0 & \text{on} \quad \partial D \times \Omega.
\end{cases} 
\end{align}
As in \cite{gunzburger2014stochastic,babuvska2007stochastic,nobile2008sparse,nobile2008anisotropic}, the following assumption is made. 
\begin{assumption}\label{inputAssumption}
The {random coefficient function} $a(\mathbf{x}, \omega)$ in system \eqref{poissonSys} has the following properties:
\begin{enumerate}
 \item There exists a positive constants $a_{\min}$ such that 
 $a_{\min} \leq a(\mathbf{x}, \omega)$,
 almost surely on $\Omega$, for all $\mathbf{x} \in D$.
 \item  $a(\mathbf{x}, \omega) = a(\mathbf{x}, \bm{\varepsilon}(\omega))$ in $\overline{D} \times \Omega$,
 where $ \bm{\varepsilon}(\omega) = (\varepsilon_1(\omega), \varepsilon_2(\omega), \ldots, \varepsilon_N(\omega))$ is a vector of real-valued uncorrelated random variables defined on a probability space $(\Omega, \mathcal{F}, \mathbb{P})$.
 \item $a(\mathbf{x}, \bm{\varepsilon}(\omega))$ is measurable with respect to $\bm{\varepsilon}$.
\end{enumerate}
\end{assumption}
For any $n=1, \ldots,N$, let $\Gamma_n := \varepsilon_n(\Omega) \subset \mathbb{R}$ and let us define the parameter space $\Gamma := \prod_{n=1}^N \Gamma_n$. 
The joint probability density function for $\{\varepsilon_n\}_{n=1}^N$ is denoted by $\rho(\bm{\varepsilon}): \Gamma \rightarrow \mathbb{R}^+$.
Several choices of coefficient functions are possible to ensure that Assumption \ref{inputAssumption} is satisfied. Some examples are given in \cite{gunzburger2014stochastic}. The choice made in this work is described in the next section.
\subsection{Karhunen-Lo\`{e}ve Expansion}\label{sectionKL}
The Karhunen-Lo\`{e}ve (KL) series expansion has been widely used in the field of uncertainty quantification to represent random fields, such as the coefficient function in system \eqref{poissonSys}, as an infinite sum of random variables \cite{li2008fourier,schevenels2004application,frauenfelder2005finite}. For numerical simulations, due to the finite computational capability available, the series is truncated, resulting in the following approximation
\begin{align}\label{truncatedKLfora}
a(\mathbf{x}, \bm{\varepsilon}(\omega)) \approx \mathbb{E}[a(\mathbf{x}, \cdot)] + \sum_{n=1}^{N} \sqrt{\lambda_n} b_n({\mathbf{x}}) \varepsilon_{n}(\omega),
\end{align}
where, for $n=1, \ldots,N$, $\lambda_n$ and $b_n$ are, respectively, the eigenvalues and eigenfunctions of the covariance function of $a(\mathbf{x}, \bm{\varepsilon}(\omega))$.
\begin{remark}\label{Gaussian_input}
We assume that the input random processes approximated with a truncated KL expansion are Gaussian, hence the random variables $\{ \varepsilon_n(\omega)\}_{n=1}^N$ are standard independent  and identically distributed. This allows us to generate realizations of the truncated KL expansion by Monte Carlo sampling from Gaussian distributions.
\end{remark}
To bound the coefficient function away from zero, the KL expansion is used on the logarithm of $a(\mathbf{x}, \bm{\varepsilon}(\omega)) - a_{\min}$ rather than on the function itself. Thus, let us define the function $\gamma(\mathbf{x}, \bm{\varepsilon}(\omega))$ as follows
\begin{align}\label{gamma}
 \gamma(\mathbf{x}, \bm{\varepsilon}(\omega)) := \log(a(\mathbf{x}, \varepsilon(\omega)) - a_{\min}).
\end{align}
The above definition implies that
\begin{align}\label{logeq}
a(\mathbf{x}, \bm{\varepsilon}(\omega)) = a_{\min} + \exp \Big( \gamma(\mathbf{x}, \bm{\varepsilon}(\omega)) \Big). 
\end{align}
Using a truncated KL expansion to approximate $\gamma(\mathbf{x}, \bm{\varepsilon}(\omega))$, we get 
\begin{align}\label{KL_gamma}
  \gamma(\mathbf{x}, \bm{\varepsilon}(\omega)) \approx \mu_{\gamma} + \sum_{n=1}^{N_{}} \sqrt{\lambda_n} b_n({\mathbf{x}}) \varepsilon_n(\omega),
\end{align}
where $\mu_{\gamma} := \mathbb{E}[\gamma(\mathbf{x}, \cdot)]$, and $N$ is the dimension of $\Gamma$.
Note that, in Eq. \eqref{KL_gamma}, ($\lambda_n$,$b_n$) are the eigenpairs associated with the covariance function of $\gamma(\mathbf{x}, \bm{\varepsilon}(\omega))$.
Moreover, according to Remark \ref{Gaussian_input}, when the truncated KL is used to approximate $\gamma(\mathbf{x}, \bm{\varepsilon}(\omega))$, the assumption of Gaussian distribution is made on $\gamma(\mathbf{x}, \bm{\varepsilon}(\omega))$, and consequently $a(\mathbf{x}, \bm{\varepsilon}(\omega))$ has a log-normal distribution, as shown in Eq. \eqref{logeq}.
It follows from Eq. \eqref{logeq} that $a(\mathbf{x}, \bm{\varepsilon}(\omega))$ is approximated by
\begin{align}\label{coeff}
a(\mathbf{x}, \bm{\varepsilon}(\omega)) \approx a_{\min} + \exp \Big(\mu_{\gamma} + \sum_{n=1}^{N} \sqrt{\lambda_n} b_n({\mathbf{x}}) \varepsilon_{n}(\omega) \Big).
\end{align}
The eigenvalues and eigenfunctions in Eq. \eqref{coeff} are obtained solving the generalized eigenvalue problem
\begin{align}\label{continuousEVP}
 \int_D C_{\gamma}(\mathbf{x},\widehat {\mathbf{x}})b_n(\mathbf{x}) d\mathbf{x} = \lambda_n b_n(\widehat{\mathbf{x}}),
 \end{align}
 where $C_{\gamma}(\mathbf{x},\widehat {\mathbf{x}})$ is the covariance function of the field $\gamma(\mathbf{x}, \bm{\varepsilon}(\omega))$.
 The covariance structure is generally unknown and usually a specific covariance function is assumed. 

In this work, the solution of Eq. \eqref{continuousEVP} is obtained according to a Galerkin approach, following the procedure presented, for instance,  in \cite{schevenels2004application,huang2001convergence}.
Assume that the physical domain $D$ is discretized with a regular finite element grid $\mathcal{T}_h$ of size $h$ \cite{ciarlet2002finite,brenner2007mathematical,aulisa2018construction}, and let $J_h$ be the total number of degrees of freedom.
If $\{\phi_j\}_{j=1}^{J_h}$ denotes the global nodal basis associated with $\mathcal{T}_h$, the eigenfunction $b_n$ of the covariance function of $\gamma(\mathbf{x}, \bm{\varepsilon}(\omega))$ in Eq. \eqref{coeff} is approximated by its nodal
interpolant
\begin{align}\label{bnInterp}
b_n(\mathbf{x}) \approx \sum_{j=1}^{J_h} b_{j,n} \phi_j(\mathbf{x}),
\end{align}
where $b_{j,n} := b_n(\mathbf{x}_j)$, with $\mathbf{x}_j$ being the $j$-th degree of freedom of $\mathcal{T}_h$. 
A substitution of Eq. \eqref{bnInterp} in Eq. \eqref{continuousEVP} transforms the continuous problem into a finite dimensional one
\begin{align}\label{almostContEVP}
  \sum_{j=1}^{J_h} b_{j,n} \Big( \int_D C_{\gamma}(\mathbf{x},\widehat {\mathbf{x}})\phi_j(\mathbf{x}) d\mathbf{x} - \lambda_n\phi_j(\widehat{\mathbf{x}}) \Big) = 0.
\end{align}
Next, Eq. \eqref{almostContEVP} is multiplied by $\phi_i$, and integrated with respect to $\widehat{\mathbf{x}}$, to obtain
\begin{align}\label{discreteEVP}
   \sum_{j=1}^{J_h} b_{j,n} \Big( \int_D \int_D C_{\gamma}(\mathbf{x},\widehat {\mathbf{x}})\phi_j(\mathbf{x}) \phi_i(\widehat{\mathbf{x}}) d\mathbf{x} d\widehat{\mathbf{x}} - \lambda_n \int_D \phi_j(\widehat{\mathbf{x}}) \phi_i(\widehat{\mathbf{x}}) d\widehat{\mathbf{x}} \Big) = 0.
\end{align}
Let us define the two real symmetric $J_h \times J_h$ matrices $\bm{C}$ and $\bm{M}$ by
\begin{align}\label{matricesEVP}
 C_{ij} := \int_D \int_D C_{\gamma}(\mathbf{x},\widehat {\mathbf{x}})\phi_j(\mathbf{x}) \phi_i(\widehat{\mathbf{x}}) d\mathbf{x} d\widehat{\mathbf{x}}, \qquad M_{ij} := \int_D \phi_j(\widehat{\mathbf{x}}) \phi_i(\widehat{\mathbf{x}}) d\widehat{\mathbf{x}}.
\end{align}
Note that $\bm{M}$ is the standard finite element mass matrix, and is positive definite.
With these matrices, Eq. \eqref{discreteEVP} can be rewritten as the vector equation
\begin{align}\label{vectorEVP}
\bm{C} \bm{b}_n = \lambda_n \bm{M} \bm{b}_n,
\end{align}
Because $\bm{C}$ is symmetric, its eigenvectors $\bm{b}_n$ associated with distinct eigenvalues are orthogonal with respect to the mass matrix $\bm{M}$, that is $\bm{b}_n^T \bm{M} \bm{b}_m = 0$ if $m \neq n$.
This orthogonality property implies that the approximations of the functions $b_n(\mathbf{x})$ in Eq. \eqref{bnInterp} are orthogonal in $L^2(D)$. Orthonormality can be obtained dividing these functions by their $L^2(D)$ norm.
Note that orthonormality in $L^2(D)$ is a requirement for the functions $b_n(\mathbf{x})$. 

For the numerical tests, the generalized eigenvalue problem in Eq. \eqref{vectorEVP} is implemented in the in-house finite element code FEMuS \cite{femus-web-page}, and solved with the SLEPc library \cite{hernandez2005slepc}.
\subsection{Quantities of Interest}\label{sectionQoI}
The aim of the paper is to approximate probability density functions of functionals of the {random PDE solution}, which we refer to as quantities of interest. 
Let $u_{J_h}(\mathbf{x}, \bm{y})$ be the approximate solution of system \eqref{poissonSys}, obtained with a {SFE method}.
The methods used to solve such a system are discussed in the next section.
Given $u_{J_h}(\mathbf{x}, \bm{y})$, examples of {random} quantities of interest include the spatial average over the physical domain
 \begin{align}\label{QoI2}
 \mathcal{Q}_u(\bm{\varepsilon}) = \dfrac{1}{|D|} \int_D u_{J_h}(\mathbf{x}, \bm{\varepsilon}) d\mathbf{x},
 \end{align}
 the integral of the square
  \begin{align}\label{QoI3}
 \mathcal{Q}_u(\bm{\varepsilon}) = \int_D \Big(u_{J_h}(\mathbf{x}, \bm{\varepsilon})\Big)^2 d\mathbf{x},
 \end{align}
 or the maximum over the physical domain
 \begin{align}\label{QoI1}
 \mathcal {Q}_u(\bm{\varepsilon}) = \max\limits_{\mathbf{x} \in D} u_{J_h}(\mathbf{x}, \bm{\varepsilon}).
 \end{align}
 Note that the dependence on the physical variable $\mathbf{x}$ is eliminated in the quantities of interest, as $\mathcal{Q}_u(\bm{\varepsilon})$ is a scalar {random variable} that only depends on the stochastic variable $\bm{\varepsilon}$.
\section{Solution of the {random PDE}}\label{SFEMs}
The {PDE with random parameters} \eqref{poissonSys} is solved numerically using stochastic finite element methods.
A detailed description of {SFE methods for PDEs} with random input data such as the one in system \eqref{poissonSys} can be found in the review article \cite{gunzburger2014stochastic} and references therein. 
The {SFE methods} used here are described in the next two sections, together with the procedures to compute stochastic quantities such as moments, once the approximate {SFE} solution is obtained.
\subsection{The Monte Carlo Method}\label{MCMsection}
The first {SFE method} considered is the classical Monte Carlo {(MC)} method, which is a stochastic sampling method \cite{gunzburger2014stochastic}.
With MC, $M$ points $\{\bm{\varepsilon}_m\}_{m=1}^M$ are chosen randomly in the parameter domain $\Gamma$, and system \eqref{poissonSys} is solved independently for each of these points.
In such a way, $M$ realization of the solution of the {random PDE} are obtained and $M$ uncoupled finite element systems are solved. 
For more details on MC and error estimates, see \cite{fishman2013monte,metropolis1949monte,gunzburger2014stochastic}. From now on, the number $M$ will refer to the number of MC samples.
 According to Remark \ref{Gaussian_input}, the MC samples are randomly drawn from Gaussian distributions. The mean $\mu_{\gamma}$ and standard deviation $\sigma_{\gamma}$ of $\gamma(\mathbf{x}, \bm{\varepsilon}(\omega))$ in Eq. \eqref{gamma} are not in general zero and one, respectively. Hence, if one chooses samples for the KL expansion from a standard Gaussian distribution, because of the presence of $\mu_{\gamma}$ in the KL and $\sigma_{\gamma}$ in the covariance function, the result is the same as sampling from a non-standard Gaussian distribution with mean $\mu_{\gamma}$  and standard deviation $\sigma_{\gamma}$.
 If $\bm{\varepsilon}_m$ denotes one of the $M$ samples obtained from the Monte Carlo Method,
 let $u_{J_h}(\mathbf{x}, \bm{\varepsilon}_m)$ be the MC solution of system \eqref{poissonSys} associated with the sample $\bm{\varepsilon}_m$ and let $\mathcal{Q}_u(\bm{\varepsilon}_m)$ be a value of a quantity of interest obtained from the realization $u_{J_h}(\mathbf{x}, \bm{\varepsilon}_m)$.
 Then $\mathcal{Q}_u(\bm{\varepsilon}_m)$ is not a {random variable} but rather just a scalar quantity.
 On the other hand, the quantity of interest $\mathcal{Q}_u(\bm{\varepsilon})$ is a function of a random variable, and a random variable itself.
 Hence stochastic quantities such as moments and cumulants can be computed. In the framework of MC, the mean of  $\mathcal{Q}_u(\bm{\varepsilon})$ is approximated with Monte Carlo integration as follows
 \begin{align}\label{meanApprox}
  \mathbb{E}\Big[\mathcal{Q}_u(\bm{\varepsilon})\Big] \approx \dfrac{1}{M}\sum_{m=1}^M \mathcal{Q}_u(\bm{\varepsilon}_m) := \mu_{\mathcal{Q}_u}.
  \end{align}
For $l > 1$, the $l$-th moment is approximated in the same fashion by
\begin{align}\label{momentApprox}
  \mathbb{E}\Big[\Big(\mathcal{Q}_u(\bm{\varepsilon}) \Big)^l\Big] \approx \dfrac{1}{M} \sum_{m=1}^M \Big(\mathcal{Q}_u(\bm{\varepsilon}_m) \Big)^l.
  \end{align}
  An estimate of the variance is given by
  $$
  \mathbb{E}\Big[\Big(\mathcal{Q}_u(\bm{\varepsilon}) - \mathbb{E}\Big[\mathcal{Q}_u(\bm{\varepsilon})\Big] \Big)^2\Big] \approx \dfrac{1}{M} \sum_{m=1}^M \Big(\mathcal{Q}_u(\bm{\varepsilon}_m) - \mu_{\mathcal{Q}_u}\Big)^2.
  $$
  The accuracy of Monte Carlo integration improves with the number of samples $M$, but the sampling error increases with the magnitude of the input variance \cite{gunzburger2014stochastic}. For this reason, here the MC is used only for small values of the input variance and to validate the stochastic Galerkin method, which is later employed for larger values of the input variance.

\subsection{The Stochastic Galerkin Method}
The second {SFE method} used to obtain an approximation of the solution of system \eqref{poissonSys} is a stochastic Galerkin {(SG)} method.
 With SG, the stochastic function space is approximated using a Galerkin procedure, as it is done for the physical function space in standard finite element methods \cite{ghanem1991stochastic,babuska2004galerkin,babuvska2005solving}.
Let us assume that the parameters are independent, so that the input joint probability function $\rho(\bm{\varepsilon})$ can be expressed as $\rho(\bm{\varepsilon}) = \prod_{n=1}^{N} \rho_n(\varepsilon_n) $.
According to a Galerkin methodology, the infinite dimensional space $L^2_{\rho}(\Gamma)$ is approximated by the finite dimensional space
\begin{equation}\label{polySpace}
\begin{aligned} 
&\mathcal{P}_{\mathcal{J}(p)}(\Gamma) = \mbox{span}\Big \{ \prod \limits_{n=1}^{N} \varepsilon_n^{p_n} \Big| \bm{p} \in \mathcal{J}(p), \varepsilon_n \in \Gamma_n \Big \}, \\
& \mbox{where} \,\, \mathcal{J}(p) = \Big \{ \bm{p} \in \mathbb{N}^{N} \Big| \sum\limits_{n=1}^{N}p_n \leq p \Big \},
\end{aligned}
\end{equation}
which corresponds to the total degree (TD) multivariate polynomial space from \cite{gunzburger2014stochastic}. Gaussian PDFs are assumed for the input variables, hence the multivariate probabilist Hermite polynomials are employed as an orthogonal basis for $\mathcal{P}_{\mathcal{J}(p)}(\Gamma)$. 
Moreover, it holds that $\Gamma_n = \mathbb{R}$ for all $n=1,\ldots,N$, and so $\Gamma = \mathbb{R}^N$.
The multivariate polynomials are obtained in a tensor product fashion from the univariate probabilist Hermite polynomials, which are appropriately scaled so that they form an orthonormal basis with respect to the PDF
\begin{align}\label{univ_pdf}
 \rho_n(\varepsilon_n) = s(\varepsilon_n) := \dfrac{\exp{(-\varepsilon_n^2/2)}}{\sqrt{2 \, \pi}},
\end{align}
{which is a standard Gaussian}.
 The scaling is carried out by dividing the $p_n$-th univariate probabilist Hermite polynomial $H_{e_{p_n}}$ by $\sqrt{n!}$, to obtain 
\begin{align}\label{orthonorm}
\int_{\mathbb{R}} \dfrac{H_{e_{p_n}}}{\sqrt{n!}} \dfrac{H_{e_{p_m}}}{\sqrt{m!}} \rho_n(\varepsilon_n) d\varepsilon_n = \delta_{nm},
\end{align}
where $\rho_n(\varepsilon_n)$ is as in Eq. \eqref{univ_pdf} and $\delta_{nm}$ is Kronecker's delta.
Multivariate $L^2_{\rho}(\Gamma)$-orthonormal Hermite polynomials are then defined as 
\begin{align}\label{hermite}
H_{e_{\bm{p}}}(\bm{\varepsilon}) = \prod_{n=1}^{N} H_{e_{p_n}}(\varepsilon_n).
\end{align}
The {SG} approximation of the solution of system \eqref{poissonSys} is defined as \cite{gunzburger2014stochastic}
$$ u^{gSG}_{J_h \, M_p}(\mathbf{x}, \bm{\varepsilon}) = 
\sum\limits_{\bm{p} \in \mathcal{J}(p)} u_{\bm{p}}(\mathbf{x}) H_{e_{\bm{p}}}(\bm{\varepsilon}), \qquad \mbox{where} \,\,
u_{\bm{p}}(\mathbf{x}) = \sum_{j=1}^{J_h} u_{\bm{p},j} \, \phi_j(\mathbf{x}).$$
We consider quantities of the form
\begin{align}\label{QoISGM}
\mathcal{Q}_u(\bm{\varepsilon}) &= \sum\limits_{\bm{p} \in \mathcal{J}(p)} \beta_{\bm{p}} H_{e_{\bm{p}}}(\bm{\varepsilon}).
\end{align}
When $\beta_{\bm{p}}$ is the average of $u_{\bm{p}}(\mathbf{x})$, we have
\begin{align}\label{QoI2SGM}
 \mathcal{Q}_u(\bm \varepsilon) = \sum\limits_{\bm{p} \in \mathcal{J}(p)} \Big( \dfrac{1}{|D|}\int\limits_{D} u_{\bm{p}}(\mathbf{x}) d\mathbf{x} \Big)  H_{e_{\bm{p}}}(\bm{\varepsilon}),
\end{align}
as in Eq. \eqref{QoI2}.
When $\beta_{\bm{p}}$ is the integral of the square of $u_{\bm{p}}(\mathbf{x})$, we have
\begin{align}\label{QoI3SGM}
 \mathcal{Q}_u(\bm \varepsilon) =  \sum\limits_{\bm{p} \in \mathcal{J}(p)} \Big( \int\limits_{D} u^2_{\bm{p}}(\mathbf{x}) d\mathbf{x} \Big) H_{e_{\bm{p}}}(\bm{\varepsilon}),
\end{align}
As opposed to what it is obtained with MC, the quantities of interest given by the {SG method} are now polynomial functions in the variable $\bm{\varepsilon}$.
Such a feature gives a considerable advantage in terms of accuracy, because moments can now be computed exactly with appropriate quadrature rules, rather than being approximated by Monte Carlo integration.
The $l$-th moment of the $\mathcal{Q}_u(\bm{\varepsilon})$ in Eq. \eqref{QoI2SGM} is defined as
\begin{align}\label{exactMoments}
 \mathbb{E}\Big[\Big(\mathcal{Q}_u(\bm{\varepsilon})\Big)^l\Big] := \int_{\mathbb{R}^{N}} \Big(\mathcal{Q}_u(\bm{\varepsilon})\Big)^l \rho(\bm{\varepsilon}) d\bm{\varepsilon}.
\end{align}
The multivariate integral in the above equation can be computed exactly with a multidimensional Hermite quadrature rule.
Exploiting Eq. \eqref{orthonorm}, monodimensional quadrature points can be obtained as the zeros of the probabilist Hermite polynomials. The quadrature weights are given by the following formula 
\begin{align}
 w_i = \dfrac{(n-1)!}{n \, (H_{e_{n-1}}(x_i))^2},
\end{align}
where $x_i$ is the $i$-th quadrature point and $H_{e_{n-1}}$ is the $(n-1)$-th probabilist Hermite polynomial.
By construction, the weights $w_i$ have the property that $\sum_{i}^{N_q} w_i = 1$, where $N_q$ is the total number of quadrature points employed in the numerical integration.
Multidimensional quadrature points and weights can be obtained in a tensor product fashion.
For the {SG method}, given $q \in \mathbb{N}$, the coefficient function $a(\mathbf{x}, \bm{\varepsilon}(\omega))$ in system \eqref{poissonSys} is computed as \cite{gunzburger2014stochastic}
\begin{equation}
\begin{aligned}\label{dataProjection}
 & a_{SG}(\mathbf{x}, \bm{\varepsilon}) = \sum\limits_{\bm{q} \in \mathcal{J}(q)} a_{\bm{q}}(\mathbf{x}) H_{e_{\bm{q}}}(\bm{\varepsilon}), \\
 & a_{\bm{q}}(\mathbf{x}) = \int_{\Gamma} a(\mathbf{x}, \bm{\varepsilon}(\omega)) H_{e_{\bm{q}}}(\bm{\varepsilon}) \rho(\bm{\varepsilon}) d\bm{\varepsilon},
\end{aligned}
\end{equation}
with $a(\mathbf{x}, \bm{\varepsilon}(\omega))$ as in Eq. \eqref{logeq}.
In the numerical results, with the exception of the integral for the data projection in Eq. \eqref{dataProjection}, all integrals are computed by an exact quadrature rule.
Among the integrals that are computed exactly, there are also the moments in Eq. \eqref{exactMoments}.
It is worth mentioning that, if a non-standard Gaussian distribution with mean $\mu_{\gamma}$ and standard deviation $\sigma_{\gamma}$ is considered for the input variables, the change of variable $(x - \mu_{\gamma}) / \sigma_{\gamma}$ would have to be made in all integrals over the parameter space. Although, one can avoid this change of variable by taking into account the non-standard distribution in the input data, as it has been discussed in Section \ref{MCMsection}.

\section{The Gram-Charlier and Edgeworth expansions}\label{expansions}
The goal of this paper is to obtain the PDF $f: \mathbb{R} \rightarrow \mathbb{R}^+$ of a given quantity of interest $\mathcal{Q}_u(\bm{\varepsilon})$ using the Gram-Charlier (GC) and the Edgeworth (ED) series expansions, which express the PDF in terms of a reference PDF.
The coefficients of these expansions are functions of the cumulants of $ \mathcal{Q}_u(\bm{\varepsilon})$, which can be easily obtained once the moments of $ \mathcal{Q}_u(\bm{\varepsilon})$ have been computed. In general, if we denote by $m_l$ and $\kappa_l$ the $l$-th moment and the $l$-th cumulant respectively, the first six cumulants are expressed in terms of the moments as follows \cite{kendall1943advanced,berberan2007expressing}:
\begin{equation}
 \begin{aligned}
  & \kappa_1 = m_1, \\
  & \kappa_2 = m_2 - m_1^2, \\
  & \kappa_3 = m_3 + 2 m_1^3 - 3 m_1 \, m_2,\\
  & \kappa_4 = m_4 - 6 m_1^4 + 12 m_1^2 m_2 - 3m_2^2 - 4 m_1 m_3,\\
  & \kappa_5 = m_5 - 5 m_4 m_1 - 10 m_3 m_2 + 20 m_3 m_1^2 + 30 m_2^2 m_1 - 60 m_2 m_1^3 + 24 m_1^5,\\
  & \kappa_6 = m_6 - 6 m_5 m_1 - 15 m_4 m_2 + 30 m_4 m_1^2 -10 m_3^2 + 120 m_3 m_2 m_1 \\
  & \quad \quad -120 m_3 m_1^3 + 30 m_2^3 - 270 m_2^2 m_1^2 + 360 m_2 m_1^4 - 120 m_1^6.
 \end{aligned}
\end{equation}
In practical applications, the GC and ED series are truncated and only a finite number of terms is retained. The truncation results in an approximation of the PDF. The behavior of the truncated series as the number of terms is increased will be subject to a discussion in the next sections.

\subsection{Derivation of the expansions}
We choose the formulation given in \cite{brenn2017revisit} to formally describe the GC and ED series.
Let $X$ be any random variable with PDF $f_X(x)$.
The characteristic function $\Psi_X(t)$ of $X$ is the Fourier transform of $f_X(x)$ \cite{kendall1943advanced}.
The cumulant generating function $K_X(t)$ is defined as $K_X(t) := \ln(\Psi_X(t))$, and it holds that
\begin{align}\label{lnpsi}
 \ln(\Psi_X(t)) = \sum\limits_{l=1}^{\infty}\kappa_{X,l}\dfrac{(i t)^l}{l!},
\end{align}
where $i$ denotes the imaginary unit and $\kappa_{X,l}$ is the $l$-th cumulant of $X$, {i.e. the
cumulants are defined as the coefficients of the series expansion of $\ln(\psi_X(t))$ with respect to $it$.}
Let $\varphi$ denote {a random variable that is distributed according to} the standard Gaussian distribution, {then, according to the notation adopted in this section, we have $f_{\varphi}(x)=s(x)$, with $s(x)$ defined in Eq. \eqref{univ_pdf}. We also observe that in the present section it is important for the clarity of the presentation to specify the random variable with which a PDF is associated, whereas in the previous sections this need was not imperative.}
With the same notation used for $X$, the following equation can be obtained
\begin{align}\label{lnphi}
 \ln(\Psi_{\varphi}(t)) = \sum\limits_{l=1}^{\infty}\kappa_{\varphi,l}\dfrac{(i t)^l}{l!}.
\end{align}
Using the properties of exponentials and convergent series, Eq. \eqref{lnpsi} and Eq. \eqref{lnphi} can be combined to get
\begin{align}\label{detour}
 \Psi_X(t) = \exp \Big( \sum\limits_{l=1}^{\infty} (\kappa_{X,l} - \kappa_{\varphi,l})\dfrac{(i t)^l}{l!}\Big) \Psi_{\varphi}(t).
\end{align}
The GC and ED expansions follow from different manipulations of Eq. \eqref{detour}.

The GC is obtained with the use of Bell polynomials $B_l(x_1, \ldots, x_l)$ \cite{bell1927partition,mihoubi2008bell}. These polynomials have the property that
\begin{align}
 \exp \Big( \sum\limits_{l=1}^{\infty} x_l \dfrac{(i t)^l}{l!}\Big) = \sum\limits_{l=0}^{\infty} B_l(x_1, \ldots, x_l) \dfrac{(i t)^l}{l!},
\end{align}
with $B_0 := 1$. Therefore, Eq. \eqref{detour} becomes
\begin{align}\label{almostGC}
 \Psi_X(t) =  \Big(1 + \sum\limits_{l=1}^{\infty} B_l(\kappa_{X,1} - \kappa_{\varphi,1}, \ldots, \kappa_{X,l} - \kappa_{\varphi,l})\dfrac{(i t)^l}{l!} \Big) \Psi_{\varphi(t)}.
\end{align}
The inverse Fourier transform applied to Eq. \eqref{almostGC} gives the GC expansion
\begin{align}\label{GC}
 f_X(x) =  \Big(1 + \sum\limits_{l=1}^{\infty} B_l(\kappa_{X,1} - \kappa_{\varphi,1}, \ldots, \kappa_{X,l} - \kappa_{\varphi,l})\dfrac{(-1)^l D_x^{(l)}}{l!} \Big) s(x),
\end{align}
where $D_x^{(l)}$ represents the $l$-th derivative operator with respect to $x$ and 
\begin{align}\label{Gauss_derivs} 
  D^{(l)}_x \Big( s(x) \Big)= (-1)^l H_{e_l}(x) s(x).
  \end{align}
 The first six Bell polynomials are given by \cite{brenn2017revisit}
\begin{equation}\label{GC_coefficients_alpha}
 \begin{aligned}
& B_1(x_1) =  x_1,\\
& B_2(x_1,x_2) =  x_1^2 + x_2,\\
& B_3(x_1,x_2,x_3) = x_1^3 + 3 x_1 x_2 + x_3,\\
& B_4(x_1,x_2,x_3,x_4) = x_1^4 + 6 x_1^2 x_2 + 3 x_2^2 + 4 x_1 x_3 + x_4,\\
& B_5(x_1,x_2,x_3,x_4,x_5) = x_1^5 + 10 x_1^3 x_2 + 15 x_1 x_2^2 + 10 x_1^2 x_3 + 10 x_2 x_3+ 5 x_1 x_4 + x_5 ,\\
& B_6(x_1,x_2,x_3,x_4,x_5,x_6) = x_1^6 + 15 x_1^4 x_2 + 45 x_1^2 x_2^2 + 15 x_2^3 + 20 x_1^3 x_3 + 60 x_1 x_2 x_3 \\
& \qquad \qquad \qquad \qquad \qquad \quad + 10 x_3^2 + 15 x_1^2 x_4 + 15 x_2 x_4 + 6 x_1 x_5 + x_6.
 \end{aligned}
\end{equation}
Note that when $X$ is a standardized variable, $\kappa_{X,0} = 0$ and $\kappa_{X,1} = 1$.
Since $\kappa_{\varphi,0} = 0$, $\kappa_{\varphi,1} = 1$, and $\kappa_{\varphi,l} = 0$ for $l \geq 3$, Eq. \eqref{GC} reduces to
\begin{align}\label{GCstandardized}
 f_X(x) =  \Big(1 + \sum\limits_{l=3}^{\infty} B_l(0,0,\kappa_{X,3}, \ldots, \kappa_{X,l})\dfrac{(-1)^l D_x^{(l)}}{l!} \Big) s(x),
\end{align}

To derive ED, one assumes that $X$ can be written as the standardized sum \cite{petrov1995limit,petrov2012sums} 
\begin{align}\label{XED}
 X = \dfrac{1}{\sqrt{r}} \sum\limits_{i=1}^r \dfrac{Z_i - \mu}{\sigma},
\end{align}
where the random variables $Z_1, Z_2, \ldots, Z_r$ are independent and identically distributed, each with mean $\mu$, standard deviation $\sigma$, and $l$-th cumulant $\kappa_{Z,l}$. Note that $\kappa_{Z,1} = \mu$ and $\kappa_{Z,2} = \sigma^2$. 
Define $\nu_l = \kappa_{Z,l} / \sigma^l$, then \cite{brenn2017revisit}
\begin{align} 
\kappa_{X,0}=0, \qquad \kappa_{X,1}=1, \qquad \kappa_{X,l} = \dfrac{\nu_l}{r^{\frac{l}{2}-1}}, \quad l \geq 3.
\end{align}
Substituting the above equation in Eq. \eqref{detour}, we get
\begin{align}\label{almostED}
 \Psi_X(t) = \exp \Big( \sum\limits_{l=3}^{\infty} \dfrac{\nu_l}{r^{\frac{l}{2}-1}}\dfrac{(i t)^l}{l!}\Big) \Psi_{\varphi}(t).
\end{align}
With a proper shift of the summation index and an inverse Fourier transform, the above series gives the ED expansion
\begin{align}\label{ED}
 f_X(x) =  \Big(1 + \sum\limits_{l=1}^{\infty} B_l(a_1,a_2, \ldots,a_l)\dfrac{1}{r^{l/2} \,\, l!} \Big) s(x),
\end{align}
The interested reader can consult \cite{brenn2017revisit} for more details on how the above equation is obtained. The coefficients $a_l$ are given by
\begin{align}
 a_l = \dfrac{\nu_{l+2}(-1)^{l+2}D_x^{(l+2)}}{(l+1)(l+2)}.
\end{align}
For ease of notation, if we write the series in Eq. \eqref{ED} as
\begin{align}\label{edgeworth}
f_{X}(x) = \sum\limits_{l=0}^{\infty} (-1)^l \dfrac{\vartheta_l(x)}{r^{l/2}},
\end{align} 
the first six coefficient functions $\vartheta_l(x)$ are given by \cite{blinnikov1998expansions}
\begin{equation}\label{coeffEDeasy}
 \begin{aligned}
  & \vartheta_0(x) = s(x),\\ 
  & \vartheta_1(x) = \frac{1}{3!} \nu_3 D_x^{(3)}\Big(s(x)\Big),\\
  & \vartheta_2(x) = \frac{1}{4!} \nu_4 D_x^{(4)}\Big(s(x)\Big) + \frac{1}{72} \nu^2_3 D_x^{(6)}\Big(s(x)\Big),\\
  & \vartheta_3(x) = \frac{1}{5!} \nu_5 D_x^{(5)}\Big(s(x)\Big) + \frac{1}{144} \nu_3 \nu_4 D_x^{(7)}\Big(s(x)\Big) + \frac{1}{1296} \nu_3^3 D_x^{(9)}\Big(s(x)\Big),\\
  & \vartheta_4(x) = \frac{1}{6!} \nu_6 D_x^{(6)}\Big(s(x)\Big) + \Big(\frac{1}{1152} \nu_4^2 + \frac{1}{720} \nu_3 \nu_5 \Big)D_x^{(8)}\Big(s(x)\Big)\\ \
  &\qquad \quad \,\, + \frac{1}{1728} \nu_3^2 \nu_4 D_x^{(10)}\Big(s(x)\Big) + \frac{1}{31104} \nu_3^4 D_x^{(12)}\Big(s(x)\Big).
 \end{aligned}
\end{equation}

\begin{figure}[h!]
\centering
\includegraphics[scale=0.75]{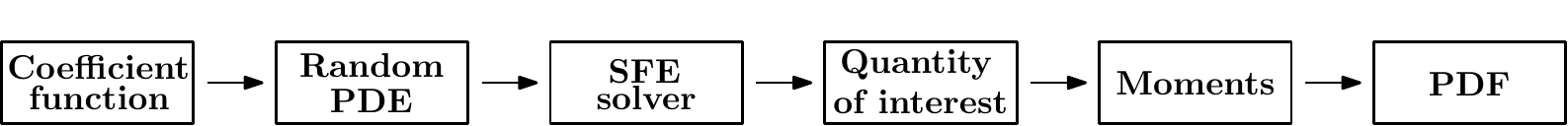}
\caption{Summary of the steps to obtain an approximation of the PDF  $f_{\mathcal{Q}_u}$.}\label{flowchart}
\end{figure}

A schematic representation of the necessary steps to get from the coefficient function $a(\mathbf{x}, \bm{\varepsilon}(\omega))$ to the approximation of the PDF $f_{\mathcal{Q}_u}$ is summarized in Figure \ref{flowchart}.
Referring to Figure \ref{flowchart}, the physical variable $\mathbf{x}$ and the stochastic variable $\bm{\varepsilon}$ are given as inputs to the {random} coefficient function $a(\mathbf{x}, \bm{\varepsilon}(\omega))$, which is used to assemble the {SFE} system associated with the {random PDE} in \eqref{poissonSys}. The system is then solved with an {SFE} method and an approximate solution is obtained. From this solution, quantities of interest are evaluated and their moments are computed, in the way associated with the {SFE} method chosen. Once the moments have been obtained, they are fed to the truncated series expansions described in this section and an approximation of the PDF $f_{\mathcal{Q}_u}$ is produced.
In the next two sections, important considerations on the convergence of the GC and ED expansions are made.
First, the concept of asymptotic expansion is discussed.

\section{Asymptotic expansions}\label{asymptexp}
Let $\{F_r(x)\}$ be a sequence of functions to be approximated by any partial sum of the series 
\begin{align}\label{wallace}
\sum\limits_{i=0}^{\infty} \dfrac{A_i(x)}{(\sqrt{r})^i}.
\end{align}
According to \cite{sedgewick2013introduction,wallace1958asymptotic}, for a given $r$, the series \eqref{wallace} is an asymptotic expansion valid to $n$ terms if
the first $n+1$ partial sums have the property
\begin{align}\label{asympt}
 \Big|F_r(x) - \sum\limits_{i=0}^{n} \dfrac{A_i(x)}{(\sqrt{r})^i}\Big| \leq \dfrac{C_n(x)}{(\sqrt{r})^{n+1}}.
\end{align}
Moreover, if $C_n(x)$ does not depend on $x$, the expansion is said to be valid uniformly in $x$.  
As pointed out by Wallace \cite{wallace1958asymptotic}, the asymptotic property is a property of finite partial sums and, for a given $r$, the series \eqref{wallace} may or may not be convergent.
If the series \eqref{wallace} converges to $F_r(x)$, we know that for every $\epsilon > 0$ there is an $R \in \mathbb{N}$ such that 
\begin{align}
  \Big|F_r(x) - \sum\limits_{i=0}^{n} \dfrac{A_i(x)}{(\sqrt{r})^i}\Big| < \epsilon, \qquad \forall n \geq R.
\end{align}
While the above inequality states that for all partial sums with $n \geq R$, the error will be less than $\epsilon$, it does not give  specific information of what happens to the error before $n$ reaches $R$. On the other hand, if $r$ is sufficiently large, an asymptotic expansion valid to $n$ terms has the property that the error will be uniformly reduced as more terms are added to the finite partial sum, from $1$ up to $n+1$ terms, after which there is no longer the guarantee that adding successive terms will provide a uniform reduction of the error \cite{sedgewick2013introduction,wallace1958asymptotic}.
With small $r$, the situation is a little different. Because the bounds $C_n(x)$ typically increase rapidly with $n$, a small value of $r$ may be unable to make the denominator grow faster than the numerator of the right hand side in Ineq. \eqref{asympt}, and only the first few terms would be improvements, as pointed out in \cite{sedgewick2013introduction,wallace1958asymptotic}. 
If a series is convergent, regardless of it being also asymptotic, the error will eventually go to zero as the number of terms in the partial sums is increased.
If a series is asymptotic but not convergent, regardless of the value of $r$, there will be a minimum error that can be achieved, which limits the accuracy of the series.
To summarize: convergent series give information on what happens to the error as the number of terms in the partial sum goes to infinity, whereas asymptotic series give information on the error as more terms are added starting from the first, and up to the $(n+1)$-th term.

{\color{black}We note that divergent asymptotic expansion have a long and useful history, especially with regards to applications. Excellent expositions on the subject are found in \cite{sedgewick2013introduction,wallace1958asymptotic}.}

\section{Considerations on convergence}\label{convconsiderations}
In \cite{kendall1943advanced} (page 152), the following sufficient condition for the convergence of the GC expansion is given
\begin{proposition}\label{GCprop}
If the PDF $f(x)$ is of bounded variation on every finite interval and
\begin{align}\label{integral}
\int_{-\infty}^{\infty} | f(x) | \exp(-x^2/4) dx < \infty,
\end{align} 
then the GC expansion of $f(x)$ is convergent.
\end{proposition}

While the integral in Proposition \ref{GCprop} is finite for PDFs with bounded support, the requirement of being of bounded variations on every finite interval may be in general harder to satisfy.
We could not find in the literature any explicit sufficient condition addressing the convergence of the Edgeworth (ED) expansion.

Although, as pointed out by several authors \cite{kendall1943advanced,cramer2016mathematical,niguez2012forecasting,fan2012probabilistic}, convergence of the GC or ED series should not be a concern, because for practical applications the series are truncated.
Of course, if the series are convergent, the error will eventually go to zero, but this is not useful for practical applications because the number of terms necessary to achieve convergence will likely be too big.
The real question is how well a truncated GC or ED series can approximate the PDF of interest.
This is directly connected to the definition of asymptotic expansion introduced in the previous section.
With an asymptotic expansion, there is the guarantee that the initial error will be reduced by adding a certain number of successive terms. This property is relevant for our purpose regardless of the convergence of the series, because only a finite number of terms after the initial approximation are considered. 
Unfortunately, the GC series is not an asymptotic expansion \cite{wallace1958asymptotic,blinnikov1998expansions,pender2014Gram,juszkiewicz1993weakly}.
On the other hand, it has been shown in \cite{cramer1928composition} that the ED expansion is asymptotic uniformly in $x$.
Even though explicit bounds on errors were not given in \cite{cramer1928composition}, the asymptotic nature of ED expansion still makes it more valuable than the GC for the applications of our interest, where the series are truncated to a small number of terms.

{\color{black}
In light of these observations, we believe that the best approach is composed of the following components. First we assume that we have in hand a stochastic Galerkin approximation of the quantity of interest. We do not include this step in the description of our algorithm because having in hand an {SG} approximation of the QoI is something one might want for any approach for estimating PDFs of QoIs.

At the heart of our algorithms is the {\em construction of a truncated Edgeworth or truncated Gram-Charlier expansion approximation of the PDF for the QoI.}
\begin{itemize}
\item[{\bf 1.}] Recursively for $l=3,\ldots$, compute the $l$-th moment of the QoI. Note that we initialize by computing three moments.
\begin{itemize}
\item Supposing we have in hand $l$ moments so obtained, we use those moments to construct an $l$-term truncated ED or GC expansion, determining cumulants from the moments.
\item The moments are obtained via exact numerical integration of the {SG} approximation of the QoI so that one need only integrate polynomials.
\end{itemize}
\end{itemize}
\noindent We need to {\em determine the ``optimal'' number of terms in the ED or GC expansions}. \begin{itemize}
\item[{\bf 2.}]This is done by comparing, each time we have incremented $l$, the $l$-term expansion with the ($l$--1)-term expansion previously determined.
\begin{itemize}
\item The comparison can be done visually or by computing, e.g., using a sampling of the two expansions, the $\ell^2$-norm of the difference between them. 
\end{itemize}
\end{itemize}
At this point we have two possibilities. 
\begin{itemize}
\item[{\bf 3a.}] {\em If the ED or GC expansion is known to be convergent or}, using the comparisons done in {\bf 2} above, {\em the expansion ``seems'' to be convergent}, we stop the recursion when the difference computed in {\bf 2} is smaller that a prescribed tolerance.
\begin{itemize}
\item In this case, the ``optimal'' number of terms in the truncated ED or GC expansion is determined by the prescribed tolerance.
\end{itemize}
\end{itemize}
\begin{itemize}
\item[{\bf 3b.}]{\em If the ED or GC expansion is known to an be asymptotic but divergent or}, using the comparisons done in {\bf 2} above, {\em the expansion ``seems'' to be divergent}, increasing the number of terms in the expansion may not result in a better approximation so that we instead determine the ``optimal'' number of terms kept in the expansion as follows.
\begin{itemize}
\item Run a Monte Carlo simulation to obtain samples of the QoI that are used to construct a histogram approximation of the PDF for the QoI, including bounds on the support of that PDF, in case it is bounded.
\begin{itemize}
\item[--] The Monte Carlo samples are determined from the SG approximation of the QoI and not by doing expensive solves of the discretized {random PDE}.
\item[--] As a result, the histogram approximation of the PDF for the QoI can be obtained at almost no cost. Even so, we do not need to take a huge number of samples because we only need to have in hand a crude histogram approximation. 
\end{itemize}
\item By comparing the crude histogram so computed with the ED or GC approximate PDFs for all values of the number of terms used, one can determine the optimal number of terms as that for which the expansion approximation is ``closest'' to the histogram.
\end{itemize}
\end{itemize}
\vskip5pt
We also remark that the number of terms kept in ED or GC expansion approximations of PDFs is often considerable lower than, e.g., kernel density estimators (KDE)  \cite{izenman1991review}. The number of terms in the latter may be as big as the number of samples taken of the {SG} approximation of the QoI. On the other hand, for ED or GC approximations of the PDF, the {SG} approximation is sampled only at the quadrature points used to approximate moment integrals, so that the number of terms kept in the expansions is not otherwise related to the number of samples. 

For the sake of simplicity, the numerical results reported in the next section do not follow every detail of the above algorithm, but those results are sufficient to illustrate the efficacy of our approach. A full implementation of the above algorithm would cast our approach in an even more favorable light as would some implementation steps that we have not discussed. An example of the latter is to take advantage of the use of {SG} approximations of the QoI to avoid approximating integrals of polynomials that are known to vanish due to the use of orthogonal polynomials in the construction of the {SG} approximation.
}

\section{Numerical Results}\label{numResults}

We consider random variables associated with the solution of \eqref{poissonSys}, obtained with the {SFE methods} discussed above.
The methods are implemented in the in-house C++ FEM solver FEMuS \cite{femus-web-page}, whereas the generalized eigenvalue problem in Eq. \eqref{vectorEVP} is solved with the SLEPc library \cite{hernandez2005slepc}.
All numerical tests consider $f \equiv - 1$ and dimension $d=2$, for simplicity. The physical domain $D$ is a unit square, with coarse grid composed of four bi-quadratic quadrilateral elements.
The mesh for the simulations is obtained by refining three times the coarse grid according to a midpoint refinement procedure.

{The random field $a(\mathbf{x}, \bm{\varepsilon}(\omega))$ is given by Eq. \eqref{coeff} with $\mu_{\gamma} = 0$ and $a_{\min} = 0.01$ and with stochastic dimension $N=2$. The covariance function of $\gamma(\mathbf{x}, \bm{\varepsilon}(\omega))$ in \eqref{logeq} is chosen to be
 \begin{align}\label{covfunc}
C_{\gamma}(\mathbf{x},\widehat {\mathbf{x}}) = \sigma_{\gamma}^2 \exp\Big[-\frac{1}{L}\Big(\sum_{i=1}^d | x_i - \widehat{x}_i|\Big)\Big],
\end{align}
where $\sigma_{\gamma}$ denotes the standard deviation of $\gamma(\mathbf{x}, \bm{\varepsilon}(\omega))$, $d$ is the dimension of the spatial variable,
and $L >0 $ is a correlation length satisfying $L \leq \mbox{diam}(D)$. The input standard deviation $\sigma_{\gamma}$ is varied during the simulations whereas the correlation length is set to $L = 0.1$. Note that although
analytic expressions are known for the eigenvalues and eigenfunctions of the covariance function Eq. \eqref{covfunc}, we do not presume such knowledge. Instead, for the numerical results given here, the eigenpairs are approximated as described in Eqs. (\ref{bnInterp} to \ref{vectorEVP}). We do so to mimic the usual situation in which eigenpairs of covariance functions are not known analytically. We also note that the approach taken in this paper also applies to anisotropic covariance functions; the choice of a correlation length $L$ that is the same for all variables $x_i$ is made merely for simplicity.}

The quantities of interest used are given in \eqref{QoI2SGM}, and \eqref{QoI3SGM}. For {the MC method}, the number of samples is $M=10^5$. Unless otherwise stated, the same amount of samples is used to obtain the crude histogram approximations. 
\Giacomo{We point out that a smaller value of $M$ can be chosen, but a fairly large number of samples allows the crude histogram to also be used for comparison with the GC and ED expansions. Within the proposed method, the histogram is obtained sampling random values of $\bm{\varepsilon}$ from a standard Gaussian distribution and plugging them in Eq. \eqref{QoISGM}, once the {SG} system has been solved. With {MC}, the histogram is obtained in the standard way.}
For {the SG method}, we choose the values of $p=4$ in Eq. \eqref{polySpace}, and $q=5$ in Eq. \eqref{dataProjection}.
\begin{remark}
When using GC, the quantity of interest is standardized before computing the moments, so the expansion adopted for the tests is Eq. \eqref{GCstandardized}.
For ED, considering Eq. \eqref{ED}, it is assumed that $r=1$ and that $Z_1 = \mathcal{Q}_u$. Hence, for both the GC and the ED  expansions, the PDF refers to the standardized quantity of interest.
\end{remark}
\subsection{Tests with nearly Gaussian output distribution}
\begin{figure}[!t]
\centering
\includegraphics[scale=0.45]{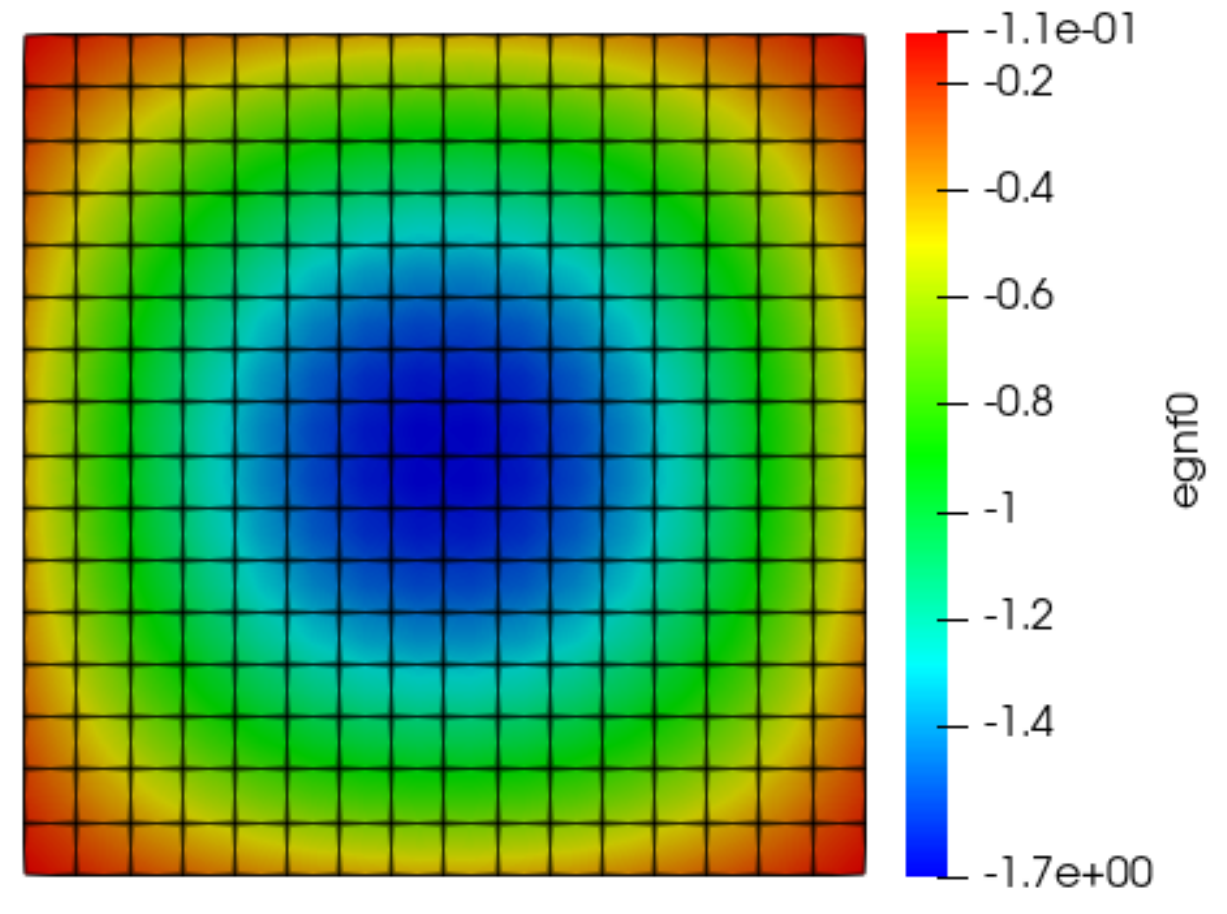}
\includegraphics[scale=0.45]{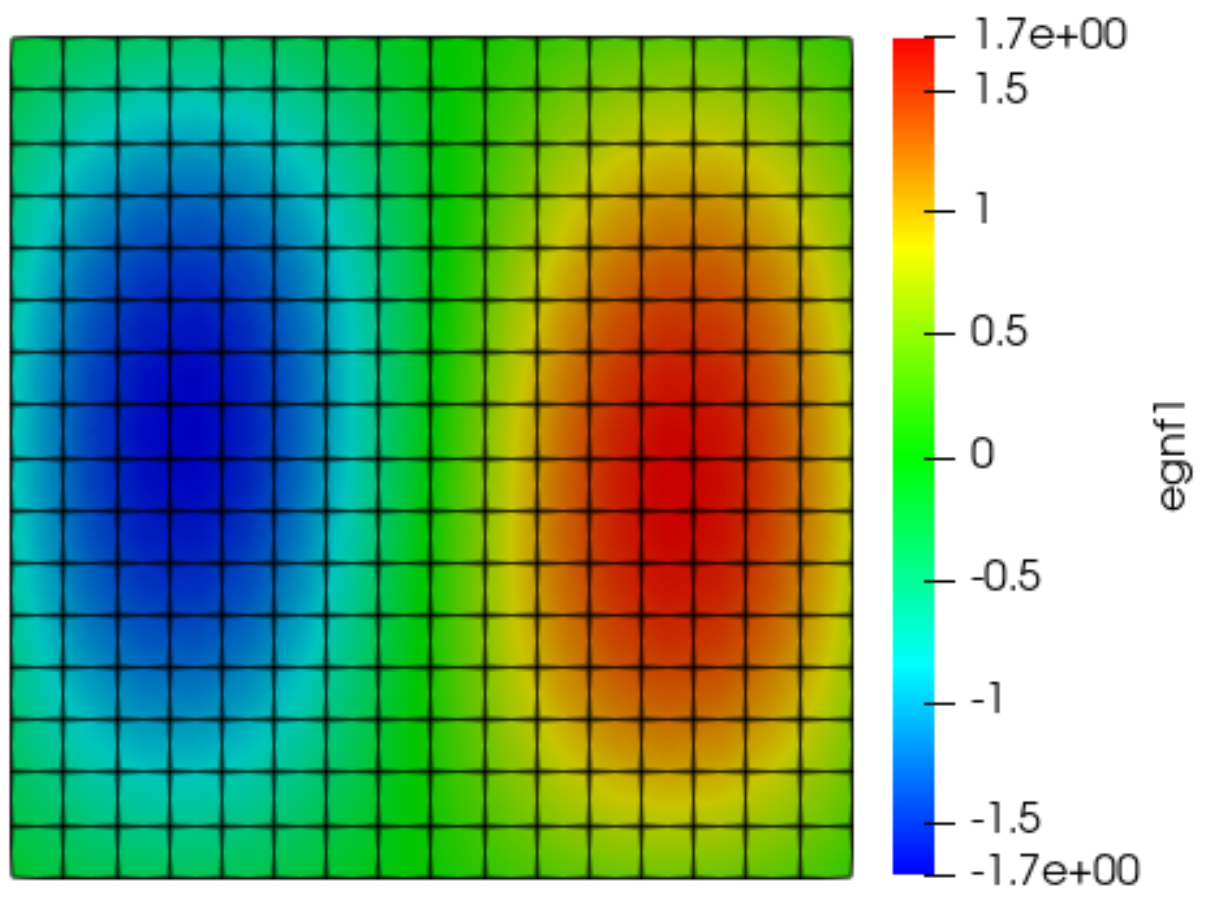}
\caption{First eigenfunction (left) and second eigenfunction (right) in the KL expansion from Eq. \eqref{coeff} for the case of $\sigma_{\gamma}=0.08$ and $N=2$.}\label{eigf}
\end{figure}

We begin with simulations where the output distribution of the quantity of interest is close to a standard Gaussian.
For the QoI in Eq. \eqref{QoI2SGM}, this is achieved with $\sigma_{\gamma} \in \{0.02, 0.04, 0.06, 0.08\}.$
Because the truncated GC and ED expansions perform successive corrections of a standard Gaussian distribution, it is expected that the expansions will perform very well in this context.
In Figure \ref{eigf}, the eigenfunctions $b_1(\mathbf{x})$ and $b_2(\mathbf{x})$ from Eq. \eqref{coeff} are reported for $\sigma_{\gamma} = 0.08$.
\begin {table}[!htb]
\setlength\tabcolsep{3.25pt} 
\begin{center}                                                               
	\begin{tabular}{|c|c|c|c|c|c|c|} \hline	
	  & \multicolumn{6}{|c|}{ QoI in Eq. \eqref{QoI2SGM} ($N=2$)} \\ \cline {2-7} 
          & \multicolumn{6}{|c|}{ $\sigma_{\gamma} = 0.08$}  \\ \cline {2-7}
  SFEM     &  $m_1$      &  $m_2$      &   $m_3$      &  $m_4$   &  $m_5$    &  $m_6$      \\ \hline 
  MC    & -3.4798e-02   & 1.2110e-03  & -4.2151e-05   & 1.4674e-06   & -5.1082e-08   & 1.7786e-09 \\ \hline 
  SG    & -3.4798e-02   & 1.2110e-03  & -4.2152e-05   & 1.4673e-06   & -5.1084e-08   & 1.7787e-09   \\ \hline 
         &  $\kappa_1$      &  $\kappa_2$     &   $\kappa_3$     &  $\kappa_4$   &  $\kappa_5$  &  $\kappa_6$      \\ \hline 
  MC    & -3.4798e-02   & 1.4246e-07  &  0   & 0   & 0   & 0        \\ \hline 
  SG    & -3.4798e-02   & 1.4356e-07  & 0   & 0  & 0   & 0     \\ \hline  
      & \multicolumn{6}{|c|}{ $\sigma_{\gamma} = 0.06$}  \\ \cline {2-7}
         &  $m_1$      &  $m_2$         &   $m_3$     &  $m_4$   &  $m_5$    &  $m_6$      \\ \hline 
  MC    & -3.4797e-02   & 1.2109e-03   & -4.2141e-05   & 1.4667e-06   & -5.1050e-08   & 1.7770e-09   \\ \hline 
  SG    & -3.4797e-02   & 1.2109e-03   & -4.2142e-05   & 1.4667e-06   & -5.1051e-08   & 1.7770e-09    \\ \hline 
         &  $\kappa_1$   &  $\kappa_2$  &   $\kappa_3$     &  $\kappa_4$   &  $\kappa_5$  &  $\kappa_6$      \\ \hline 
  MC    & -3.4797e-02   &  8.0123e-08   &   0   &  0   &  0   &  0       \\ \hline 
  SG    & -3.4797e-02   &  8.0745e-08   &   0   &  0   &  0   &  0    \\ \hline 
          & \multicolumn{6}{|c|}{ $\sigma_{\gamma} = 0.04$}  \\ \cline {2-7}
         &  $m_1$      &  $m_2$         &   $m_3$     &  $m_4$   &  $m_5$    &  $m_6$      \\ \hline 
  MC    & -3.4796e-02   & 1.2108e-03   & -4.2135e-05   & 1.4663e-06   & -5.1027e-08   & 1.7758e-09  \\ \hline 
  SG    & -3.4797e-02   & 1.2108e-03   & -4.2135e-05   & 1.4663e-06   & -5.1028e-08   & 1.7759e-09     \\ \hline 
         &  $\kappa_1$   &  $\kappa_2$  &   $\kappa_3$     &  $\kappa_4$   &  $\kappa_5$  &  $\kappa_6$      \\ \hline 
  MC   &  -3.4796e-02   &  3.5607e-08 &   0   &  0   &  0   &  0      \\ \hline 
  SG    &  -3.4796e-02   &  3.5884e-08   &   0   &  0   &  0   &  0    \\ \hline 
            & \multicolumn{6}{|c|}{ $\sigma_{\gamma} = 0.02$}  \\ \cline {2-7}
         &  $m_1$      &  $m_2$         &   $m_3$     &  $m_4$   &  $m_5$    &  $m_6$      \\ \hline 
  MC     & -3.4796e-02   & 1.2108e-03   & -4.2131e-05   & 1.4660e-06   & -5.1014e-08   & 1.7752e-09 \\ \hline 
  SG     & -3.4796e-02   & 1.2108e-03   & -4.2131e-05   & 1.4660e-06   & -5.1014e-08   & 1.7752e-09    \\ \hline 
         &  $\kappa_1$   &  $\kappa_2$  &   $\kappa_3$     &  $\kappa_4$   &  $\kappa_5$  &  $\kappa_6$      \\ \hline 
  MC    & -3.4796e-02   & 8.9013e-09   &   0   & 0   & 0   & 0        \\ \hline 
  SG    & -3.4796e-02   & 8.9705e-09    &   0   &  0   &  0   &  0    \\ \hline 
	\end{tabular}
\end{center}
\caption{Behavior of moments and cumulants of the QoI in Eq. \eqref{QoI2SGM} for $\sigma_{\gamma} \in \{0.02, 0.04, 0.06, 0.08\}$.}
\label{test1_average_N2}
\end{table}
The values of the moments and cumulants for the quantity of interest in Eq. \eqref{QoI2SGM} are shown in Table \ref{test1_average_N2}. Note that such values are associated with the non-standardized quantity of interest.
As discussed in \cite{gunzburger2014stochastic}, the sampling error of the {MC method} increases with the input standard deviation, hence we use this method only for small values of $\sigma_{\gamma}$, such as those considered in this section. The results obtained with {the MC method} are relevant because they serve as a validation of the {SG method}, that is employed in the next section for {non-Gaussian} output distributions.
Observing Table \ref{test1_average_N2}, it is fair to say that convergence has been reached for the moments, given that there is a good agreement between the {MC} and the {SG} predictions.
Therefore, the {SG method} can be considered validated by the {MC} results.
All cumulants $\kappa_l$ with $l\geq3$ are zero,
confirming that the PDF is close to a Gaussian distribution, which is the one distribution that has $\kappa_l=0$ for $l\geq3$.

\begin{figure}[!t]
  \centering
 \begin{tabular}{@{}p{0.475\linewidth}@{\quad}p{0.475\linewidth}@{}}
    \subfigimg[width=\linewidth]{I)}{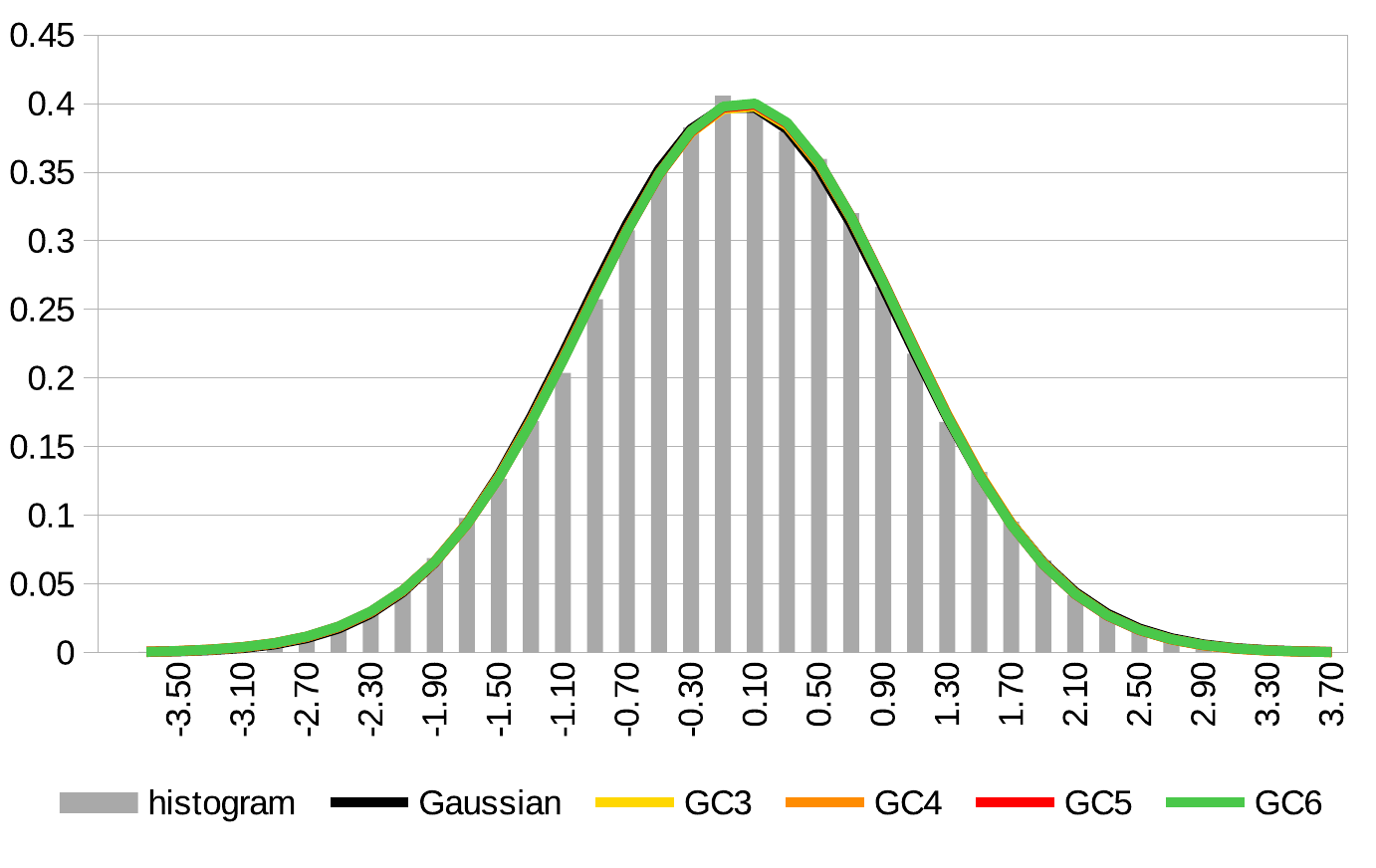}&
    \subfigimg[width=\linewidth]{II)}{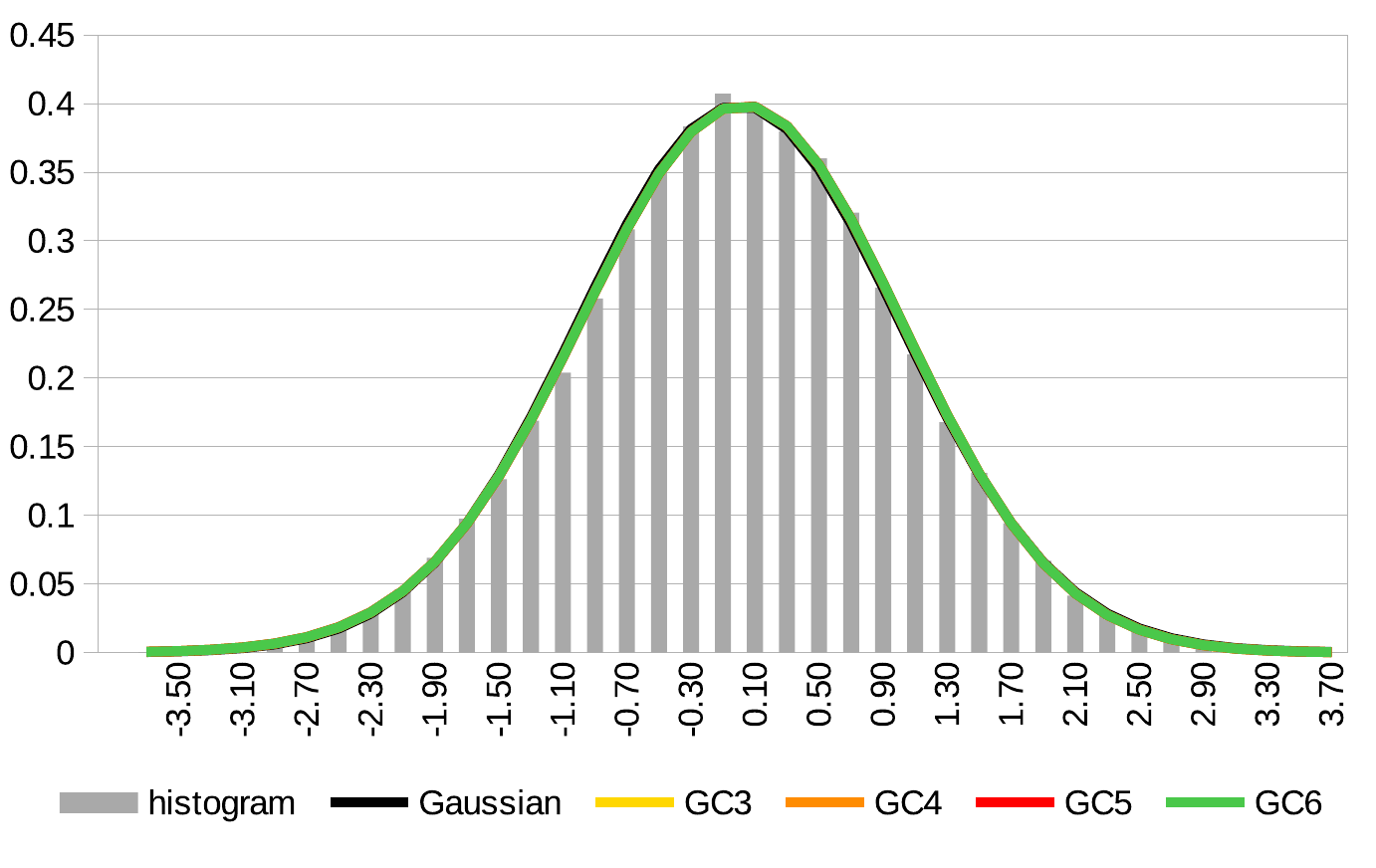}
  \end{tabular}
\caption{Gram-Charlier expansions of the standardized QoI in Eq. \eqref{QoI2SGM} for $\sigma_{\gamma}=0.08$ and $N=2$.
I): GC - MC, II): GC - SG.}
\label{test1Curves}
\end{figure}

The histograms for $\sigma_{\gamma} = 0.08$ are in Figure \ref{test1Curves}. As expected, {due to the small values of the input variance,} both suggest a nearly-Gaussian PDF:  the major
difference with an actual Gaussian distribution is that in this case the PDFs have
bounded support due to the fact that the quantity of interest is bounded. In fact,
the leading term of the expansion is a normal distribution, with the additional
terms being higher-moment contributions. If $\sigma_{\gamma}$ is small, then those contributions
are relatively small compared to the leading term, so that the overall distribution
in nearly normal.
Our tests have shown that larger values of $\sigma_{\gamma}$ cause an increased skewness on the histograms, although this effect is very weak for the values of $\sigma_{\gamma}$ considered here, and so histograms for $\sigma_{\gamma} = 0.02, 0.04$ and $0.06$ are not shown due to their strong similarity to those in Figure \ref{test1Curves}.
The increased skewness will be clearly visible when larger values of $\sigma_{\gamma}$ are considered in the next section.

In Figure \ref{test1Curves}, the truncated GC expansion is displayed for $\sigma_{\gamma}=0.08$.
GC expansions for other values of $\sigma_{\gamma}$ are not reported due to their strong resemblance to those in Figure \ref{test1Curves}. 
The notation GC$j$ means that the GC expansion in Eq. \eqref{GCstandardized} has been truncated at  $l=j$.
Truncations up to $l=6$ have been computed.
The GC expansion is a function defined on all $\mathbb{R}$, hence it could be graphed for ideally any $x$. Although, after the range of values of the quantity of interest is obtained from the crude approximation, the graph of the GC is plotted only for values of $x$ in such a range. The same is done for the truncated Edgeworth expansions.
Referring to Figure \ref{test1Curves}, we see that the GC expansions are in good agreement with the respective histogram.

\begin{figure}[!t]
  \centering
 \begin{tabular}{@{}p{0.475\linewidth}@{\quad}p{0.475\linewidth}@{}}
    \subfigimg[width=\linewidth]{I)}{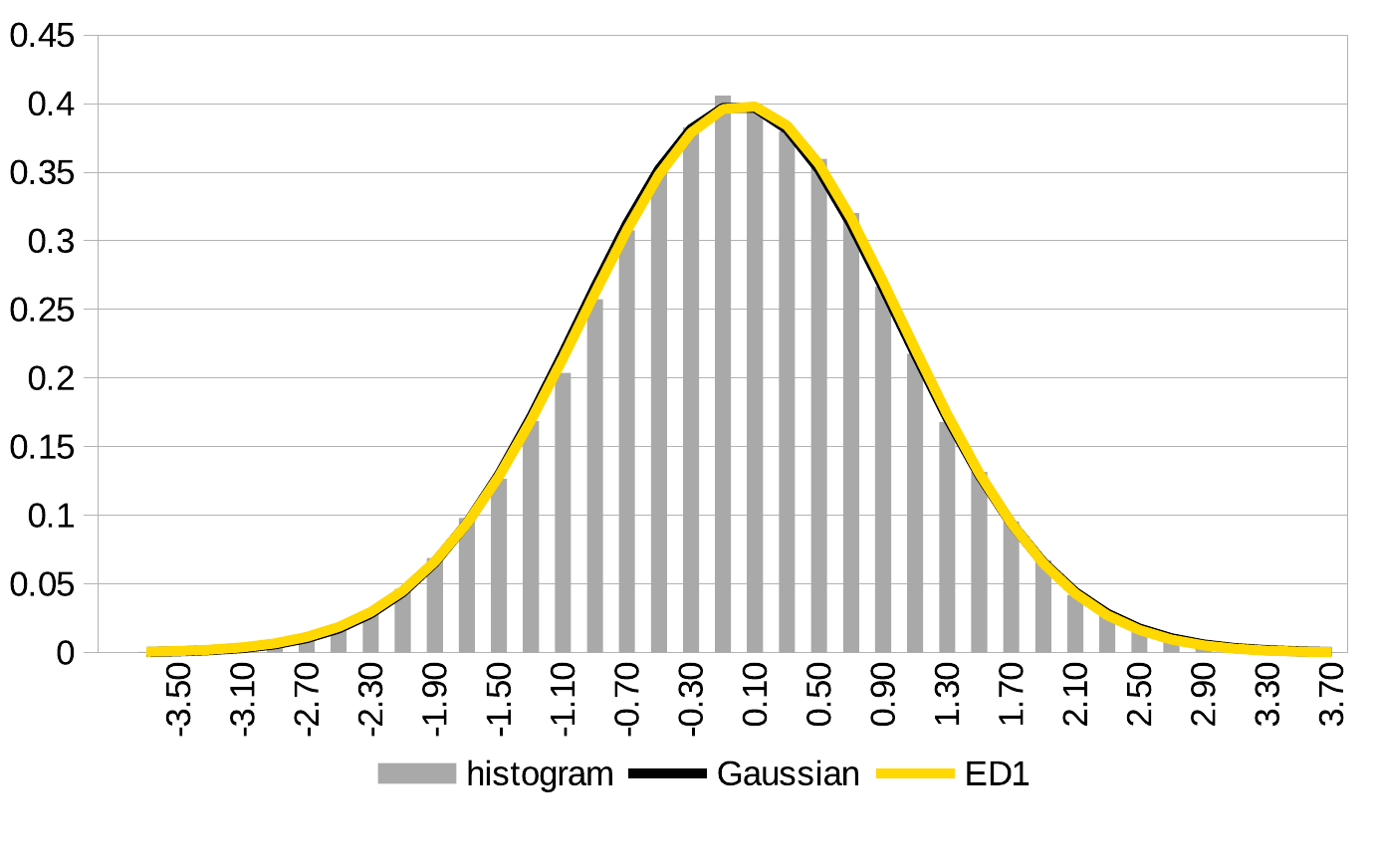} & 
     \subfigimg[width=\linewidth]{II)}{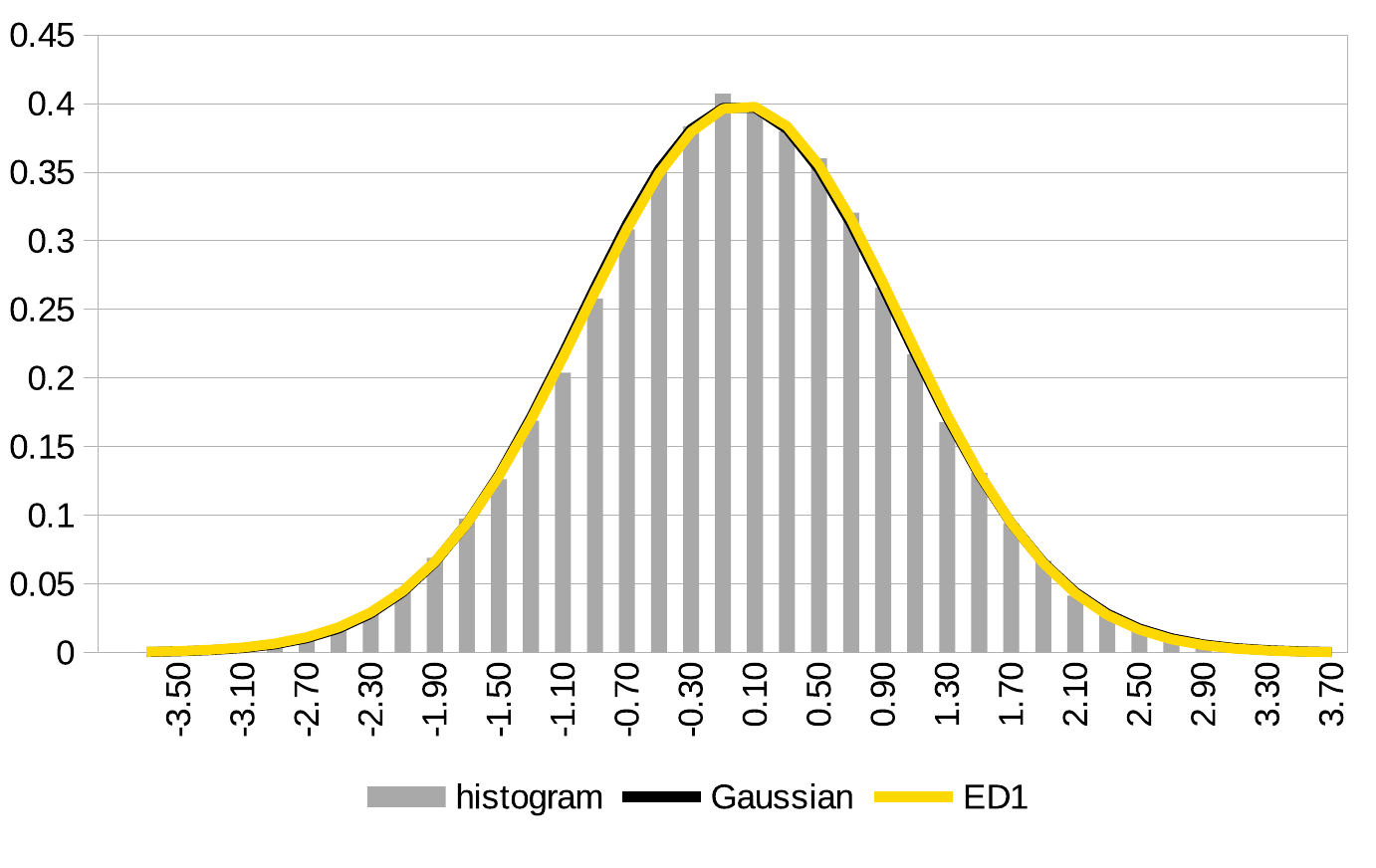}
  \end{tabular}
\caption{Edgeworth expansions of the standardized QoI in Eq. \eqref{QoI2SGM} for $N=2$ and $\sigma_{\gamma}=0.08$.
I): ED - MC, II): ED - SG.}
\label{test1_ED}
\end{figure}

In Figure \ref{test1_ED} the truncated ED expansions are shown for {MC and SG}. Only the results for $\sigma_{\gamma}=0.08$ are reported, because other values produced similar graphs.
Note that, for GC, all the truncated series computed are displayed in Figure \ref{test1Curves}, namely GC3, GC4, GC5 and GC6.
As discussed above, GC is not an asymptotic expansion and so it might be that, for instance, GC4 is a worse approximation than GC3 but GC5 is better than both GC3 and GC4.
Hence, it makes sense to display all computed curves, to have an idea of the behavior of the truncated series.
For ED, truncations up to four terms have been computed.
However, not all curves are displayed in the figures, but only those that show a monotone reduction of the error.
Such curves are typically ED1 and ED2. This scenario is expected considering the discussion in Section \ref{asymptexp}, and that $r=1$ in Eq. \eqref{XED}.
Figure \ref{test1_ED} shows great agreement of the ED expansion with the histograms, for both the {MC method and the SG method}.

When the distribution of the quantity of interest is nearly Gaussian, our results suggest that the Gram-Charlier and Edgeworth expansions can effectively be used to describe the PDF of the quantity of interest.
This is consistent with the nature of such expansions, that add successive corrections to a standard Gaussian distribution.
The natural question that arises next concerns how well GC and ED can work when the PDF to approximate is not nearly Gaussian. Numerical results are presented in the next section to address this uncertainty.
\subsection{Tests with {non-Gaussian} output distribution}
\begin {table}[!t]
\setlength\tabcolsep{3.25pt} 
\begin{center}                                                               
	\begin{tabular}{|c|c|c|c|c|c|c|c|c|} \hline	   
	  & \multicolumn{8}{|c|}{ Standardized QoI in Eq. \eqref{QoI2SGM} (SG, $N=2$)  } \\ \cline {2-9}
	  &  \multicolumn{4}{|c|}{ Moments } &  \multicolumn{4}{|c|}{ Cumulants } \\ \cline {2-9}  
 $\sigma_{\gamma}$  &  $m_{3_{STD}}$   &  $m_{4_{STD}}$   &   $m_{5_{STD}}$  &  $m_{6_{STD}}$  &  $\kappa_{3_{STD}}$  &  $\kappa_{4_{STD}}$ &   $\kappa_{5_{STD}}$ &  $\kappa_{6_{STD}}$     \\ \hline 
    1.6  & -0.793290 & 4.16802   &  -10.4848  & 46.2597  & -0.79329  & 1.16802  & -2.55194   & 7.44632 \\ \hline   
    1.4 & -0.68550   & 3.87021   & -8.49255   & 36.8592  & -0.68550  & 0.87021  & -1.63752   & 4.10686 \\ \hline
    1.2 & -0.58120   & 3.62412   & -6.80469   & 29.8402  & -0.58120  & 0.62412  & -0.99265   & 2.10040 \\ \hline
    1   & -0.47986   & 3.42448   & -5.35446   & 24.6366  & -0.47986  & 0.42448  & -0.55579   & 0.96666 \\ \hline
  	\end{tabular}
\end{center}
\caption{Moments and cumulants of the standardized QoI in Eq. \eqref{QoI2SGM} with $N=2$, for different values of the input standard deviation.}
\label{test2_spatial_average}
\end{table}
\begin {table}[!t]
\setlength\tabcolsep{3.25pt} 
\begin{center}                                                               
	\begin{tabular}{|c|c|c|c|c|c|c|c|c|} \hline	  
	  & \multicolumn{8}{|c|}{ Standardized QoI in Eq. \eqref{QoI3SGM} (SG, $N=2$)  } \\ \cline {2-9}
	  &  \multicolumn{4}{|c|}{ Moments } &  \multicolumn{4}{|c|}{ Cumulants } \\ \cline {2-9}  
 $\sigma_{\gamma}$  &  $m_{3_{STD}}$   &  $m_{4_{STD}}$   &   $m_{5_{STD}}$  &  $m_{6_{STD}}$  &  $\kappa_{3_{STD}}$  &  $\kappa_{4_{STD}}$ &   $\kappa_{5_{STD}}$ &  $\kappa_{6_{STD}}$     \\ \hline 
    3   & 0.82751 & 4.31992 & 11.3622 & 51.2865 & 0.82751 & 1.31992 & 3.08712 & 9.63985 \\ \hline
    2.9 & 0.69515 & 3.89341 & 8.62298 & 37.3751 & 0.69515 & 0.89341 & 1.67138 & 4.14151  \\ \hline
    2.8 & 0.58557 & 3.62103 & 6.81178 & 29.6887 & 0.58557 & 0.62103 & 0.95604 & 1.94421   \\ \hline
    2.7 & 0.53517 & 3.51856 & 6.08253 & 27.0037 & 0.53517 & 0.51856 & 0.73080 & 1.36114   \\ \hline
  	\end{tabular}
\end{center}
\caption{Moments and cumulants of the standardized QoI in Eq. \eqref{QoI3SGM} with $N=2$, for different values of the input standard deviation.}
\label{test2_integral}
\end{table}
%
%
%
\begin{figure}[!t]
  \centering
 \begin{tabular}{@{}p{0.475\linewidth}@{\quad}p{0.475\linewidth}@{}}
    \subfigimg[width=\linewidth]{I)}{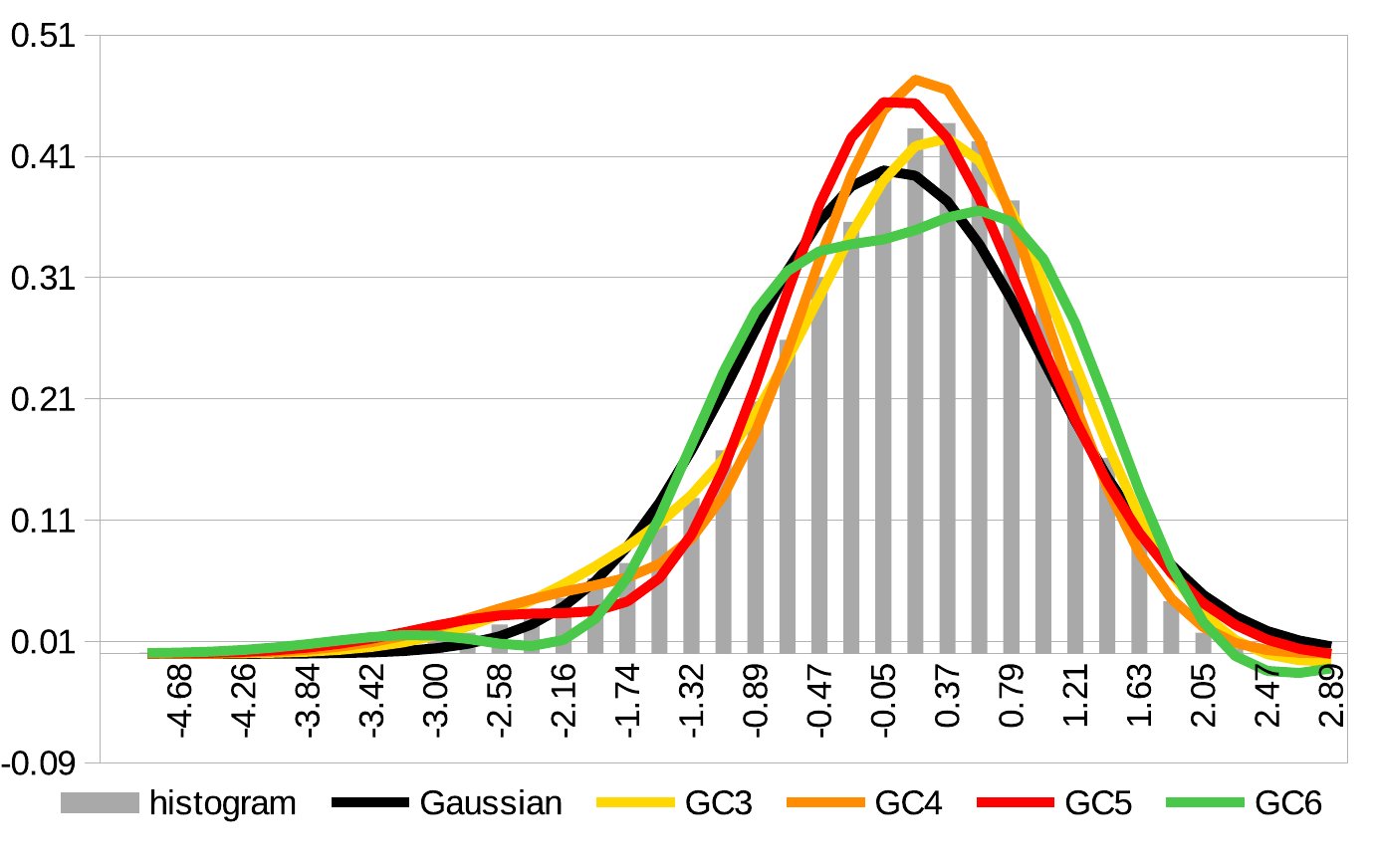}&
    \subfigimg[width=\linewidth]{II)}{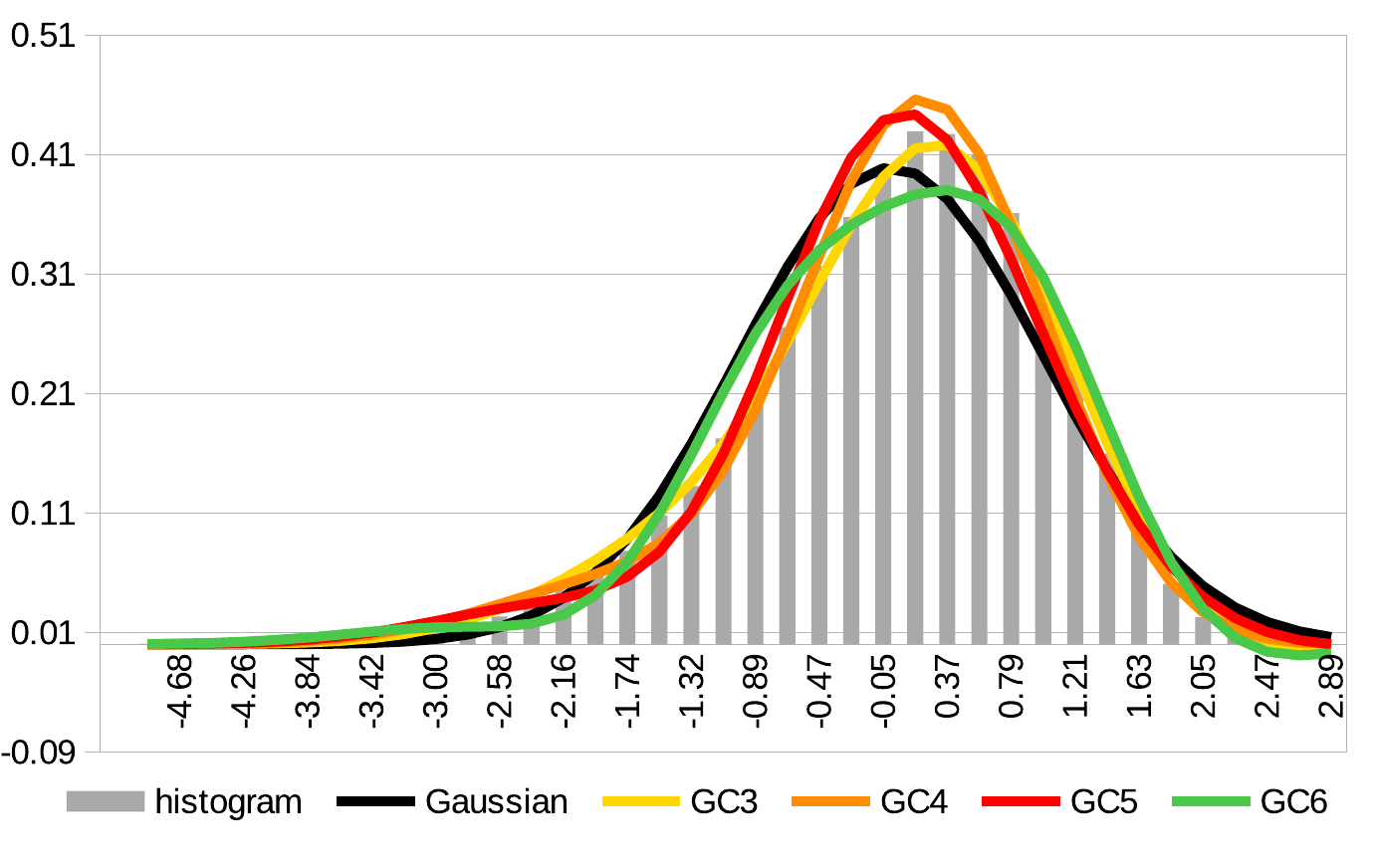} \\
    \subfigimg[width=\linewidth]{III)}{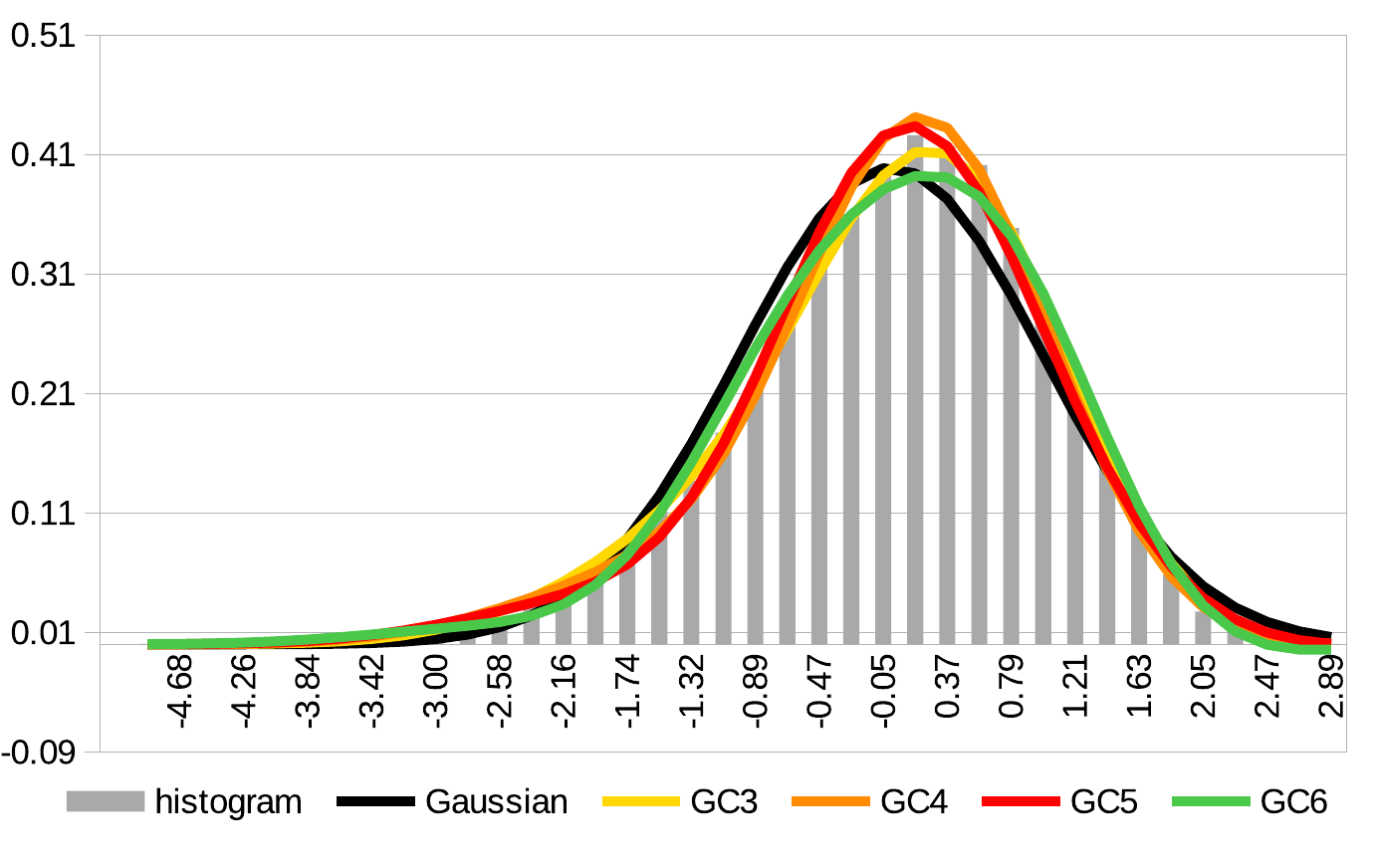} &
    \subfigimg[width=\linewidth]{IV)}{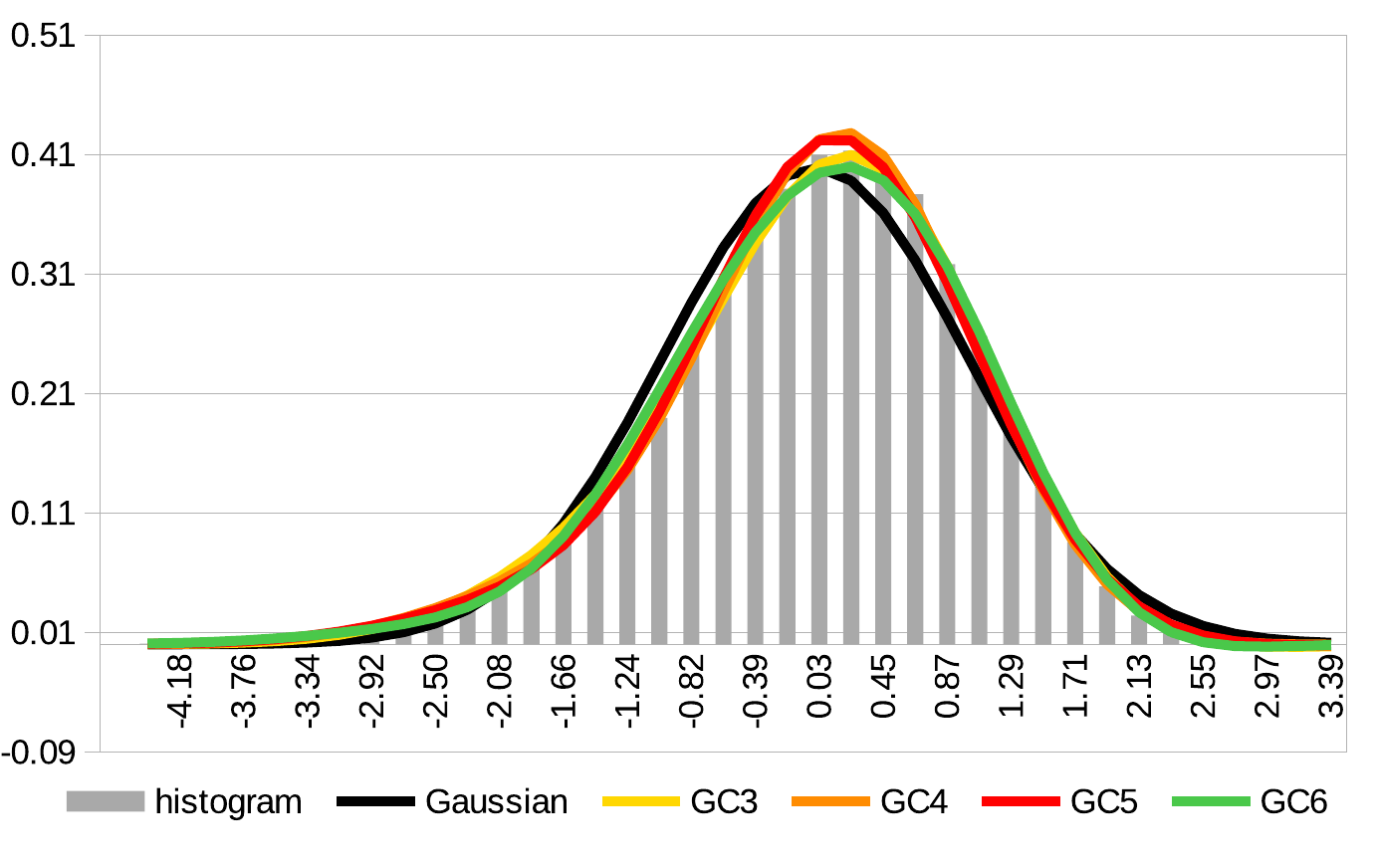} 
    \end{tabular}
\caption{Gram-Charlier expansion of the standardized QoI in Eq. \eqref{QoI2SGM} for $N=2$, obtained with the SG method.
I): $\sigma_{\gamma}=1.6$, II): $\sigma_{\gamma}=1.4$,  III): $\sigma_{\gamma}=1.2$, IV): $\sigma_{\gamma}=1$.}
\label{test2_GC_average}
\end{figure}
\begin{figure}[!t]
  \centering
 \begin{tabular}{@{}p{0.475\linewidth}@{\quad}p{0.475\linewidth}@{}}
    \subfigimg[width=\linewidth]{I)}{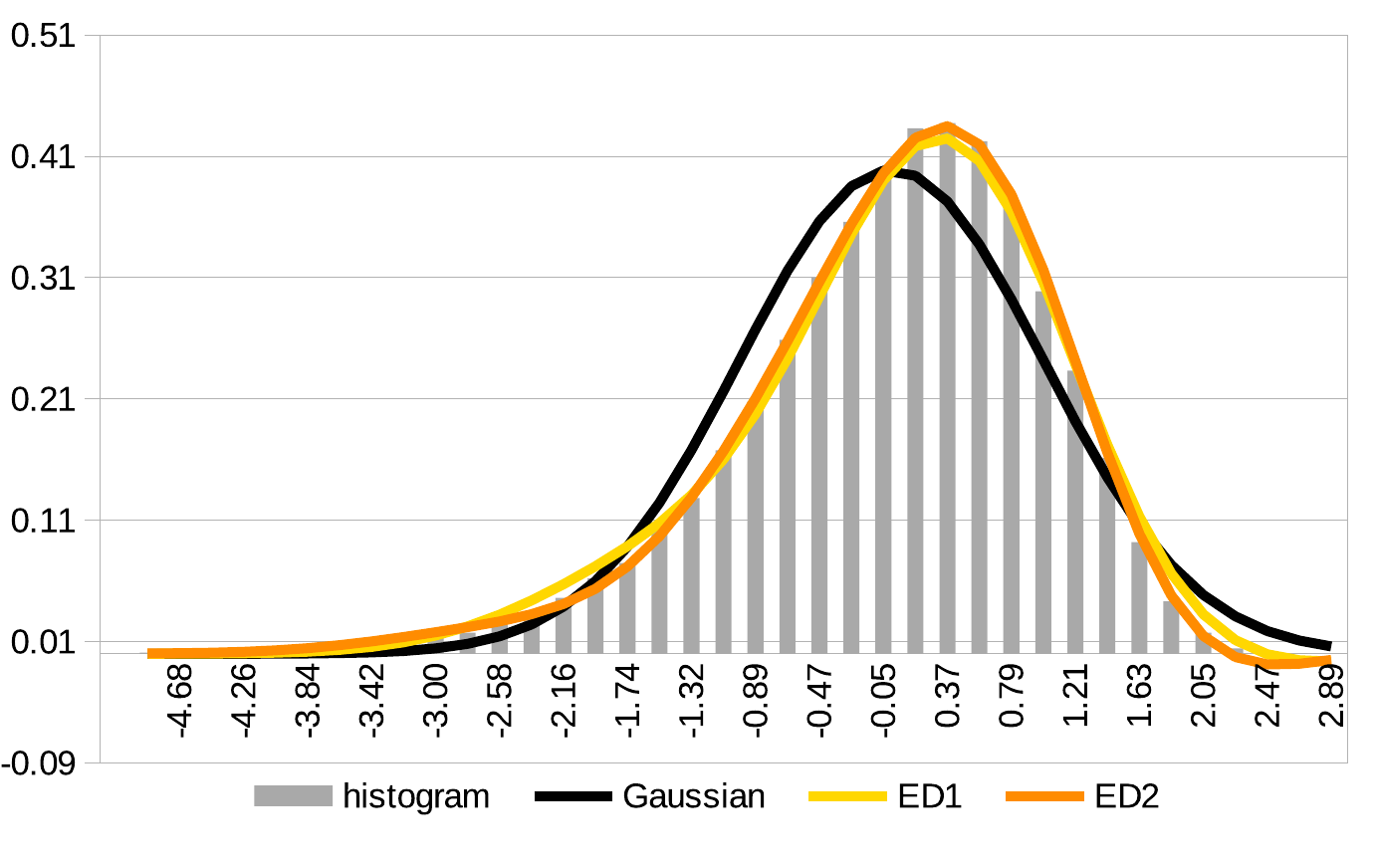}&
    \subfigimg[width=\linewidth]{II)}{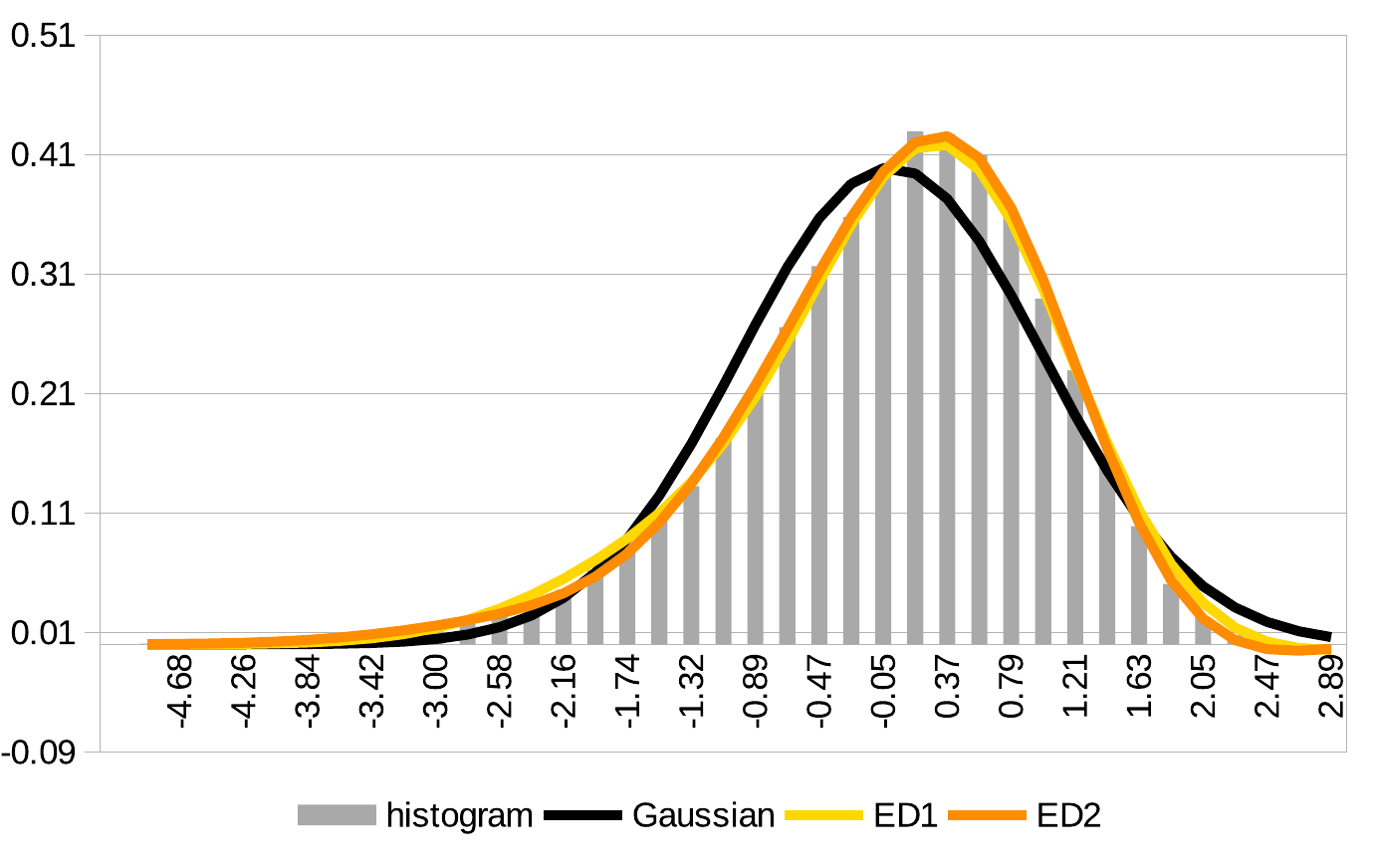} \\
    \subfigimg[width=\linewidth]{III)}{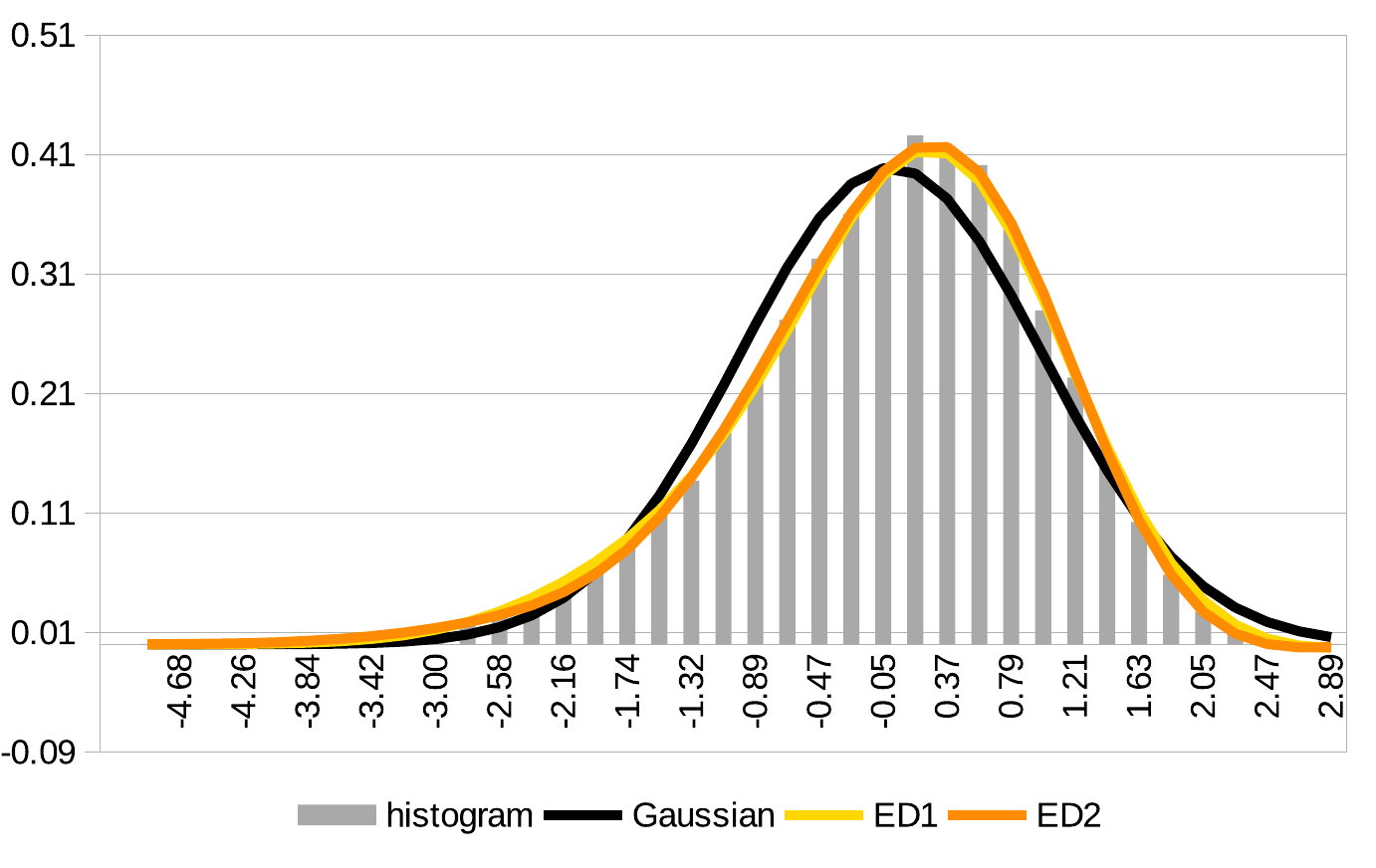} &
    \subfigimg[width=\linewidth]{IV)}{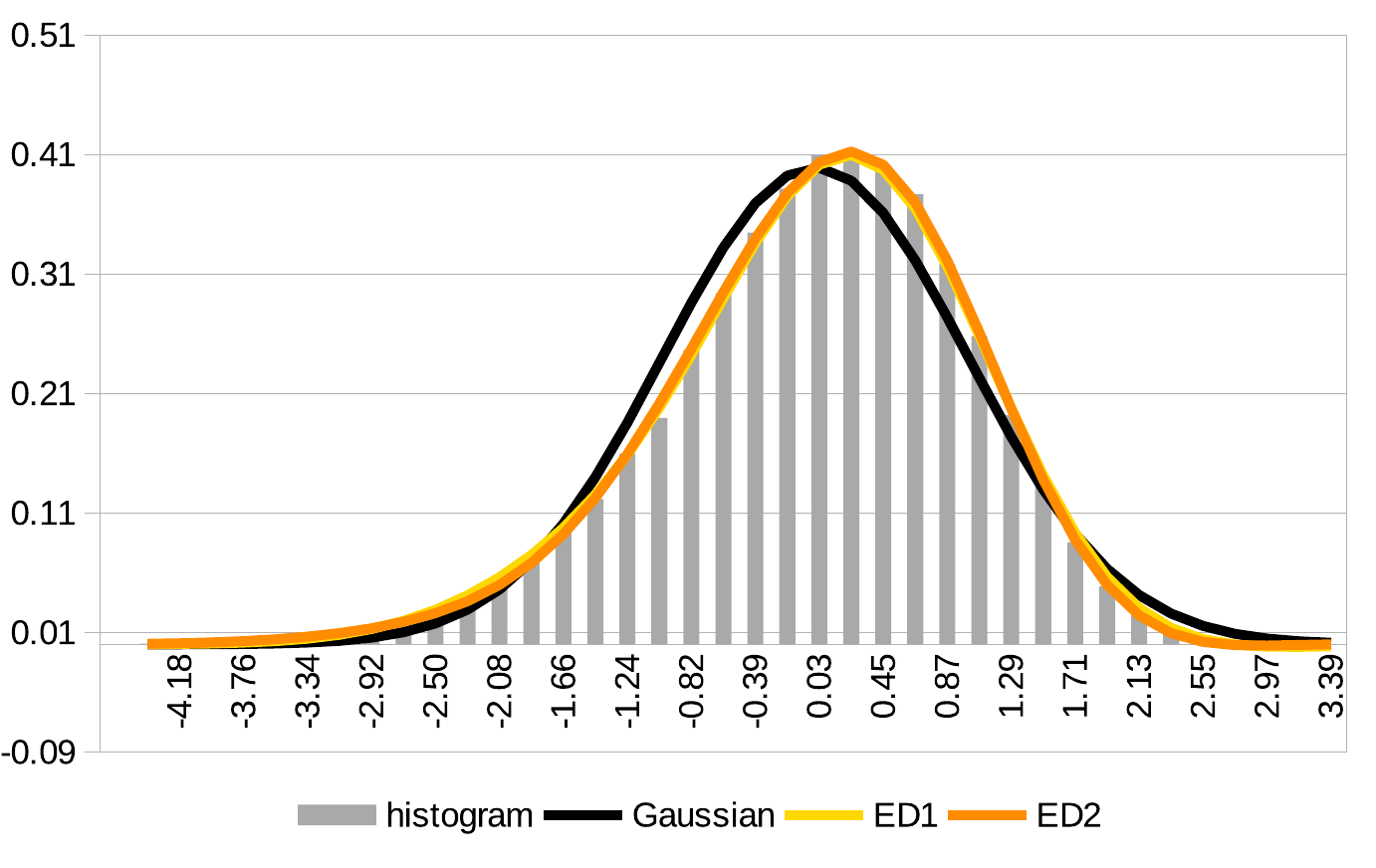} 
    \end{tabular}
\caption{Edgeworth expansion of the standardized QoI in Eq. \eqref{QoI2SGM} for $N=2$, obtained with the SG method.
I): $\sigma_{\gamma}=1.6$, II): $\sigma_{\gamma}=1.4$, III): $\sigma_{\gamma}=1.2$, IV): $\sigma_{\gamma}=1$.}
\label{test2_ED_average}
\end{figure}
We now investigate how well the GC and ED expansions can approximate a PDF that is not nearly Gaussian, {in particular we are interested in PDFs for which not just the first and the second moment are non-zero}. This is achieved with larger values of the input standard deviation than those considered in the previous section.
The QoI in \eqref{QoI2SGM} and QoI in \eqref{QoI3SGM} are considered.
With the former, $\sigma_{\gamma} \in \{1, 1.2, 1.4, 1.6 \}$, while with the latter $\sigma_{\gamma} \in \{2.7, 2.8, 2.9, 3\}$.
Only the {SG method} is used, because for the values of $\sigma_{\gamma}$ considered here, the sampling error of the {MC method} would produce results that are not accurate enough. The value $N=2$ is chosen for the stochastic dimension.
In Table \ref{test2_spatial_average} and Table \ref{test2_integral}, the values of the moments and cumulants are reported for the standardized QoIs. As the input standard deviation grows, all moments increase in magnitude, suggesting that larger values of $\sigma_{\gamma}$ give increasingly non-Gaussian distributions.
In Figure \ref{test2_GC_average}, the histograms for the QoI in \eqref{QoI2SGM} are reported: 
larger values of $\sigma_{\gamma}$ produce an increasing negative skewness.
The GC curves are in Figure \ref{test2_GC_average}.
For all values of $\sigma_{\gamma}$, GC3 provides the best approximation of the PDF, among the GC curves that have been computed. The quality of the approximation slightly decreases as the input standard deviation grows but remains overall satisfactory, especially for lower values of $\sigma_{\gamma}$.
GC4 and GC5 are also good approximations for $\sigma_{\gamma} = 1$ and $\sigma_{\gamma} = 1.2$, whereas GC6 may be considered good enough only for $\sigma_{\gamma}=1.$
In general, since the GC and ED series keep on correcting a standard Gaussian, it could be that more terms are required for a good approximation when the PDF is far from a Gaussian. 

The ED results are shown in Figure \ref{test2_ED_average}. The best approximation among the curves computed is given by ED2 which well fits all histograms. The quality of the approximation slightly decreases as $\sigma_{\gamma}$ grows, but the magnitude of this deterioration is much smaller than in the GC case.
The ED expansions are qualitatively superior to the GC curves, because ED2 is a better approximation than GC3.
Note that GC3 is by definition the same curve as ED1.

\begin{figure}[!t]
  \centering
 \begin{tabular}{@{}p{0.475\linewidth}@{\quad}p{0.475\linewidth}@{}}
    \subfigimg[width=\linewidth]{I)}{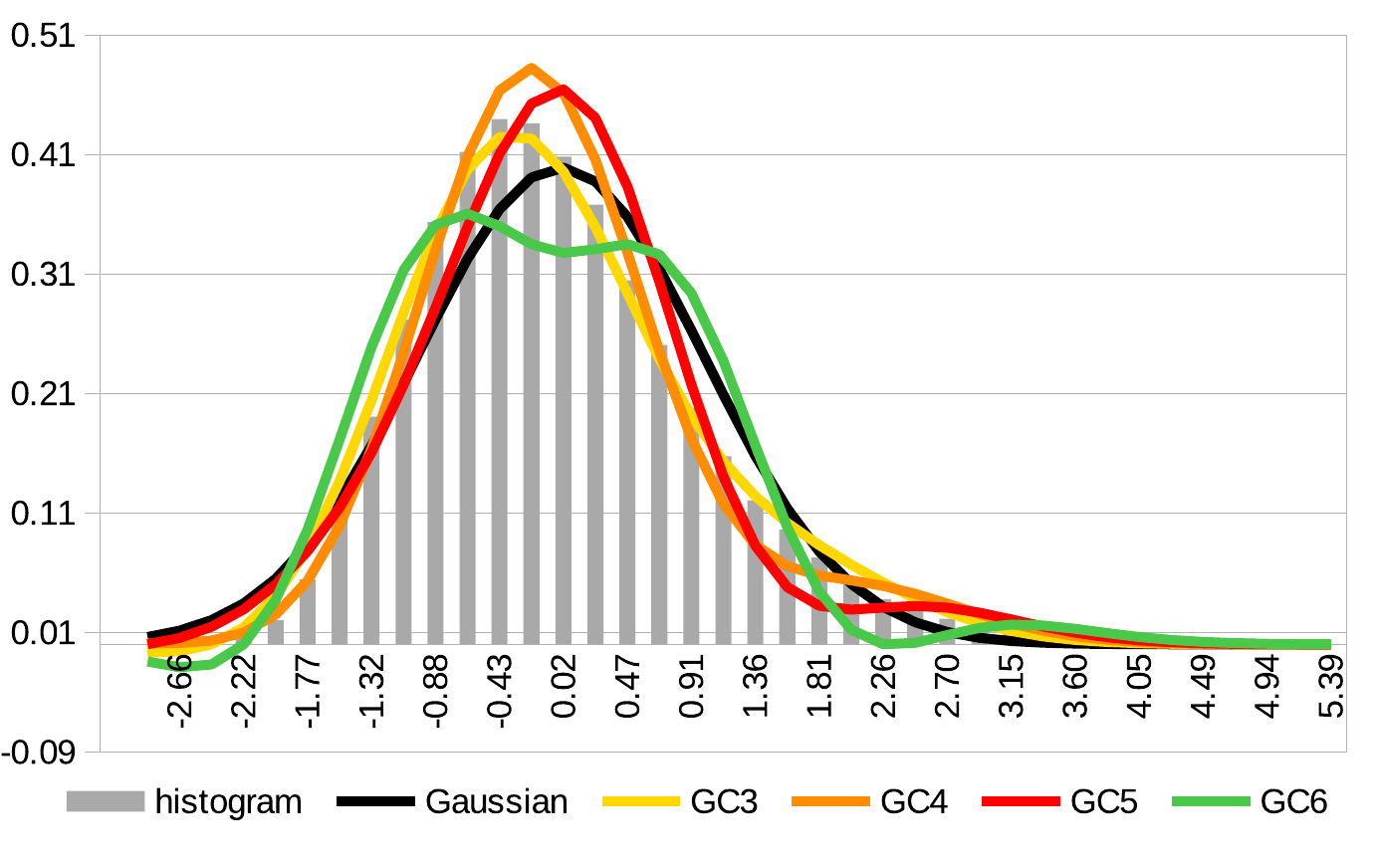}&
    \subfigimg[width=\linewidth]{II)}{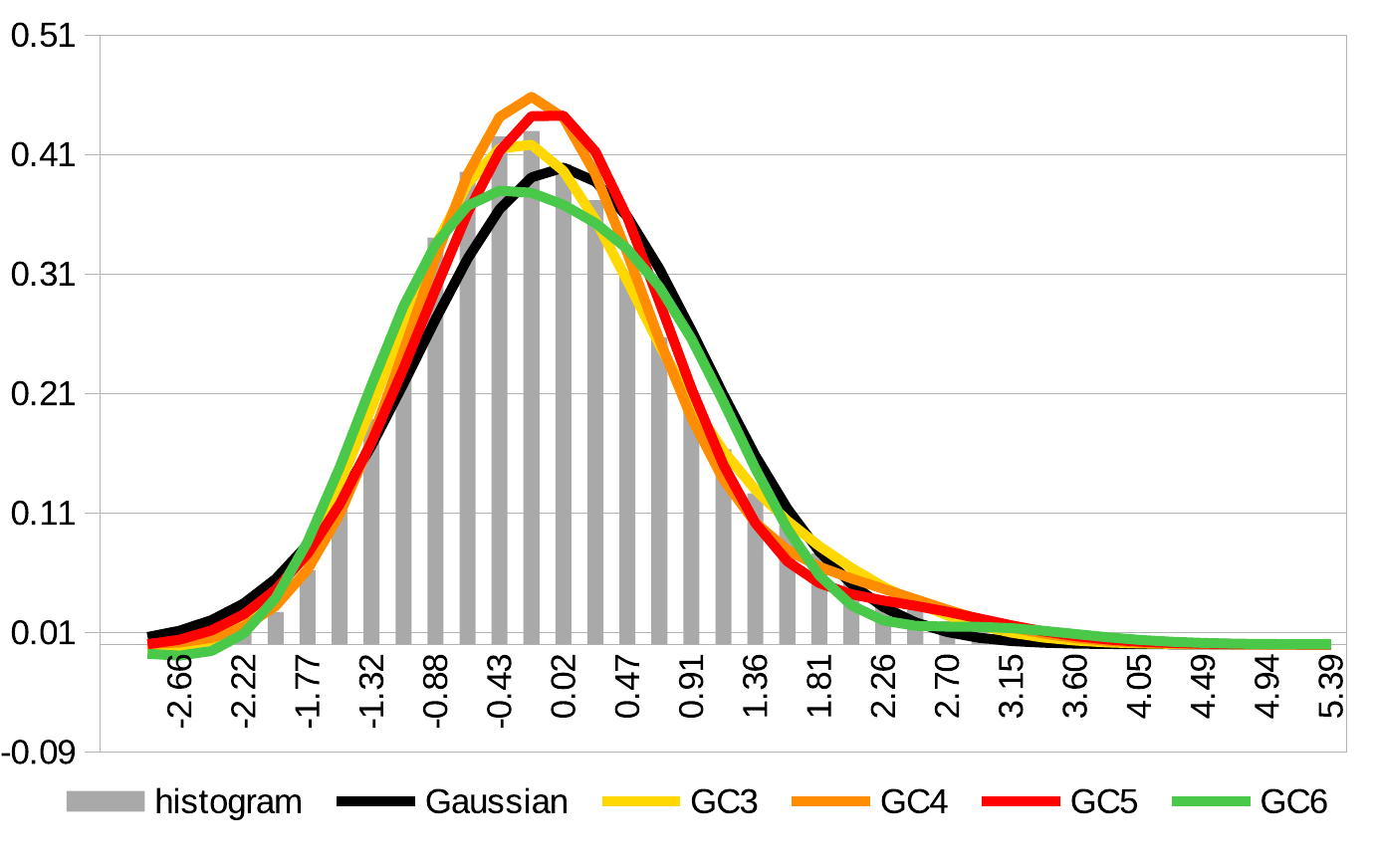} \\
    \subfigimg[width=\linewidth]{III)}{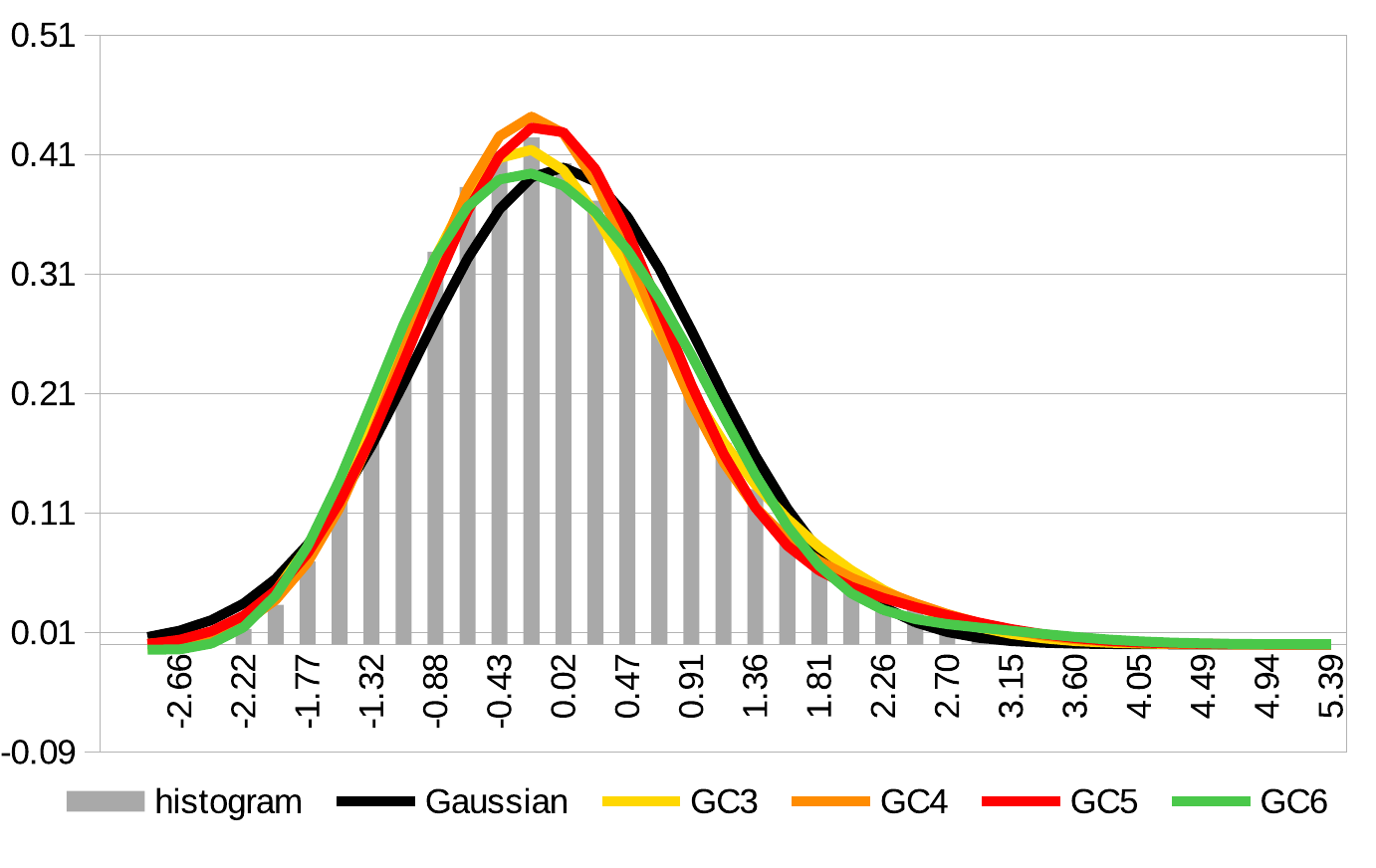} &
    \subfigimg[width=\linewidth]{IV)}{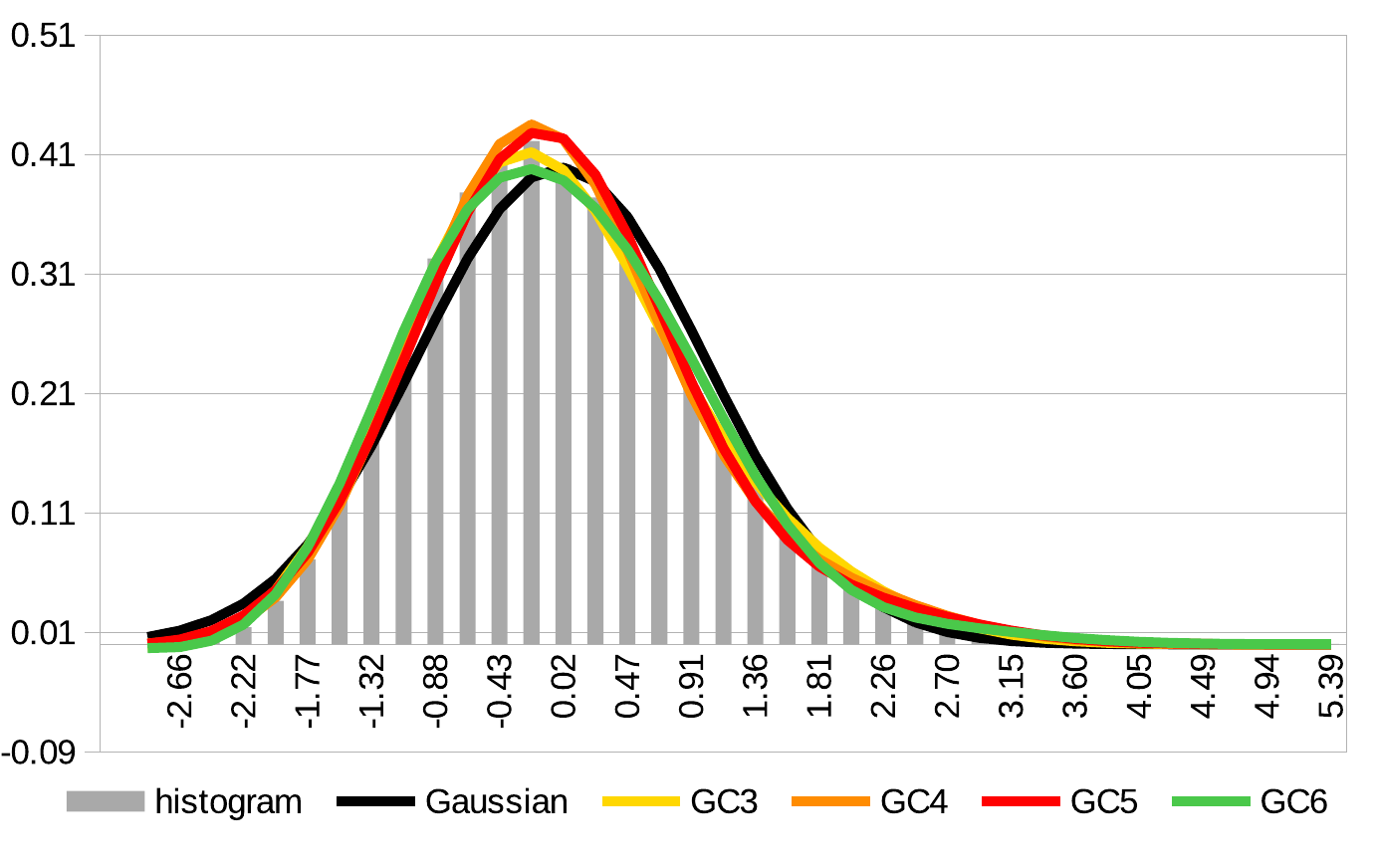} 
    \end{tabular}
\caption{Gram-Charlier expansion of the standardized QoI in Eq. \eqref{QoI3SGM} for $N=2$, obtained with the SG method.
I): $\sigma_{\gamma}=3$, II): $\sigma_{\gamma}=2.9$, III): $\sigma_{\gamma}=2.8$, IV): $\sigma_{\gamma}=2.7$.}
\label{test2_GC_integral}
\end{figure}

Next, we consider the quantity of interest given in Eq. \eqref{QoI3SGM}. The histograms are in Figure \ref{test2_GC_integral}.
Greater values of $\sigma_{\gamma}$ cause an increasing positive skewness, especially going from 
$\sigma_{\gamma}=2.9$ to  $\sigma_{\gamma}=3$. 
The GC curves are in Figure \ref{test2_GC_integral}. For $\sigma_{\gamma}=2.7$, all computed GC curves well approximate the PDF, with GC4 and GC5 lying on top of each other. 
For all the other values of $\sigma_{\gamma}$, GC3 is again the best approximation, as the other computed curves progressively deviated from the histogram. For $\sigma_{\gamma} =2.8$, GC4 and GC5 are still acceptable approximations, however they are not accurate enough for the subsequent values of $\sigma_{\gamma}$. GC6 is acceptable only for $\sigma_{\gamma}=2.7$.

\begin{figure}[!t]
  \centering
 \begin{tabular}{@{}p{0.475\linewidth}@{\quad}p{0.475\linewidth}@{}}
    \subfigimg[width=\linewidth]{I)}{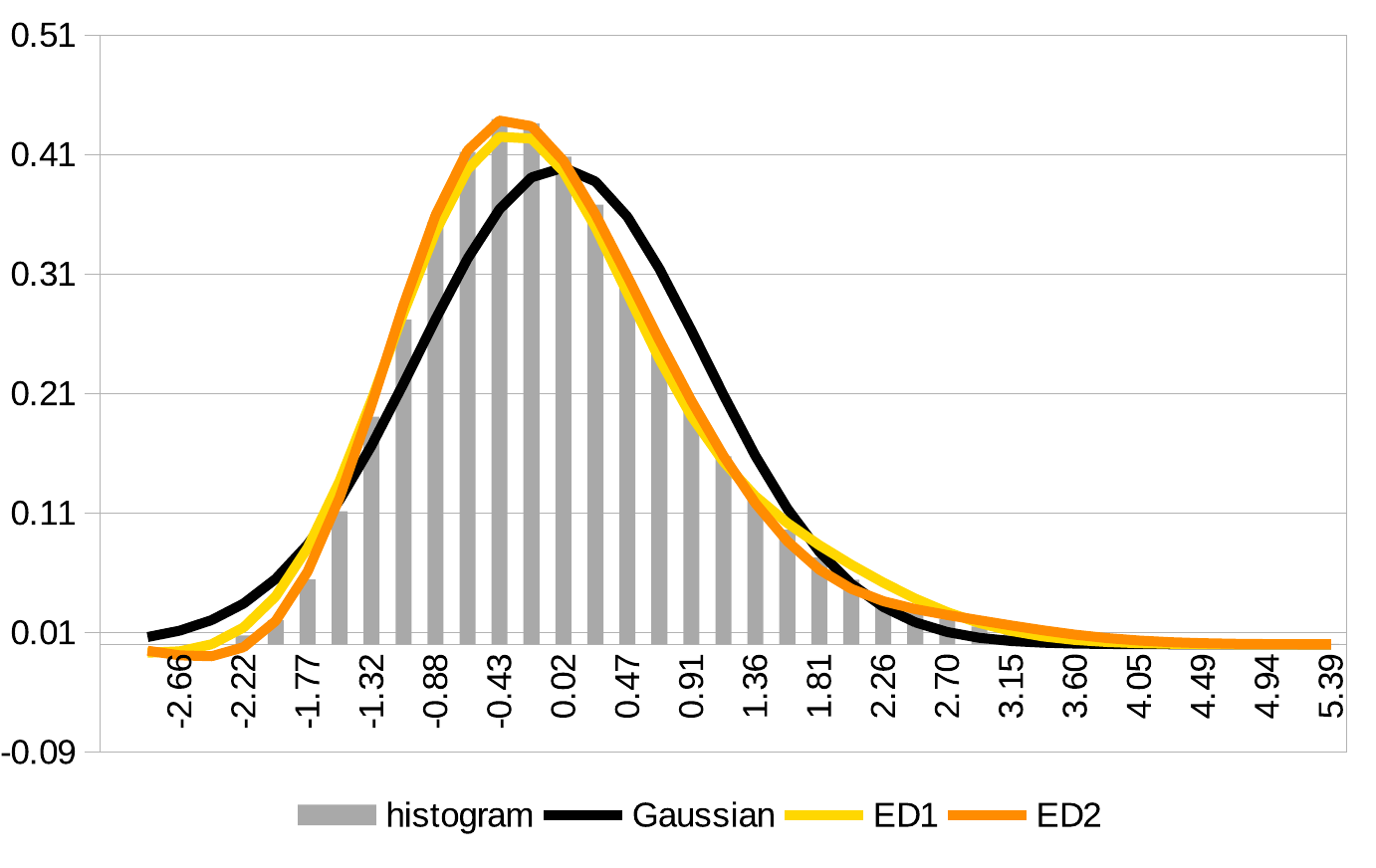}&
    \subfigimg[width=\linewidth]{II)}{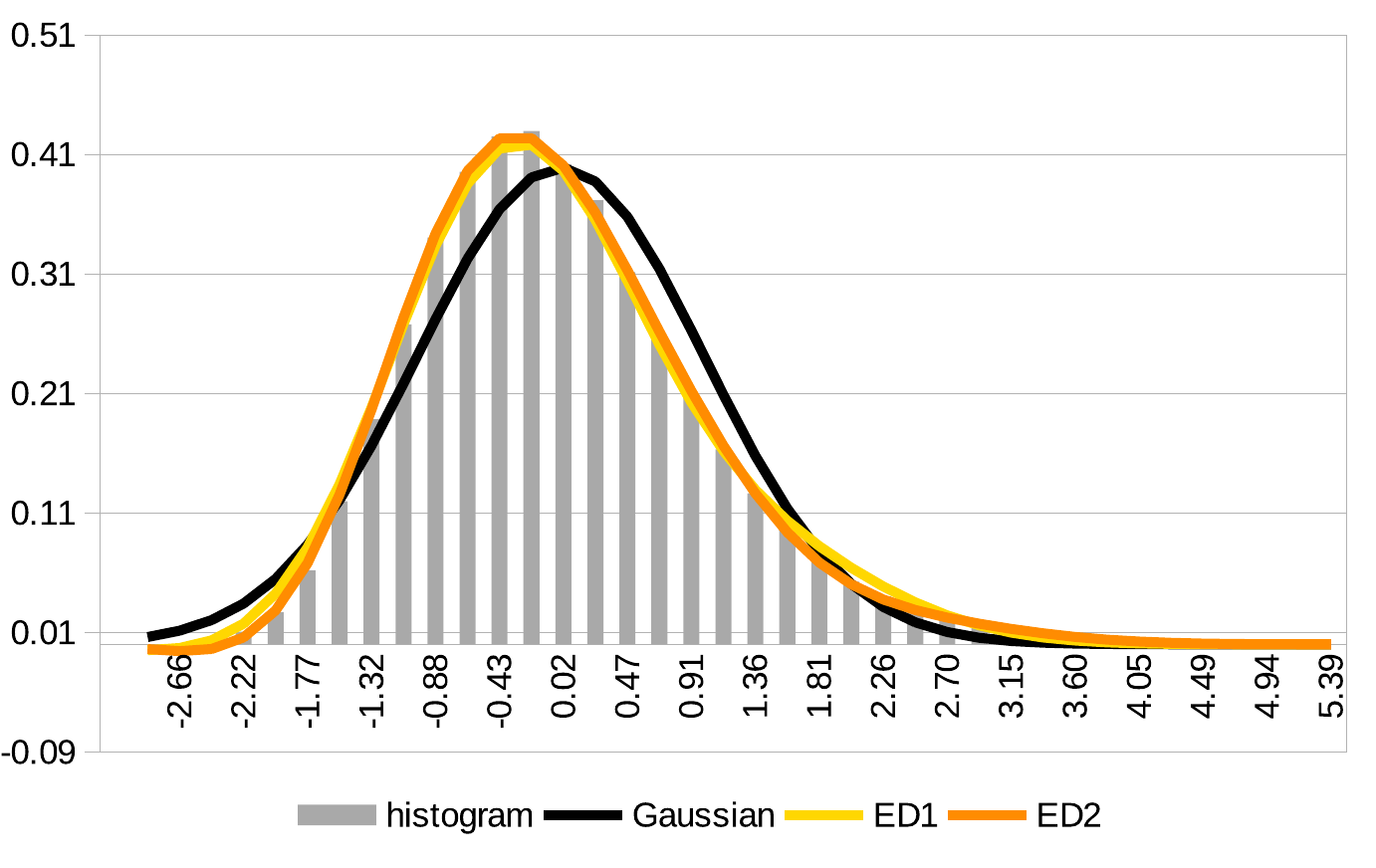} \\
    \subfigimg[width=\linewidth]{III)}{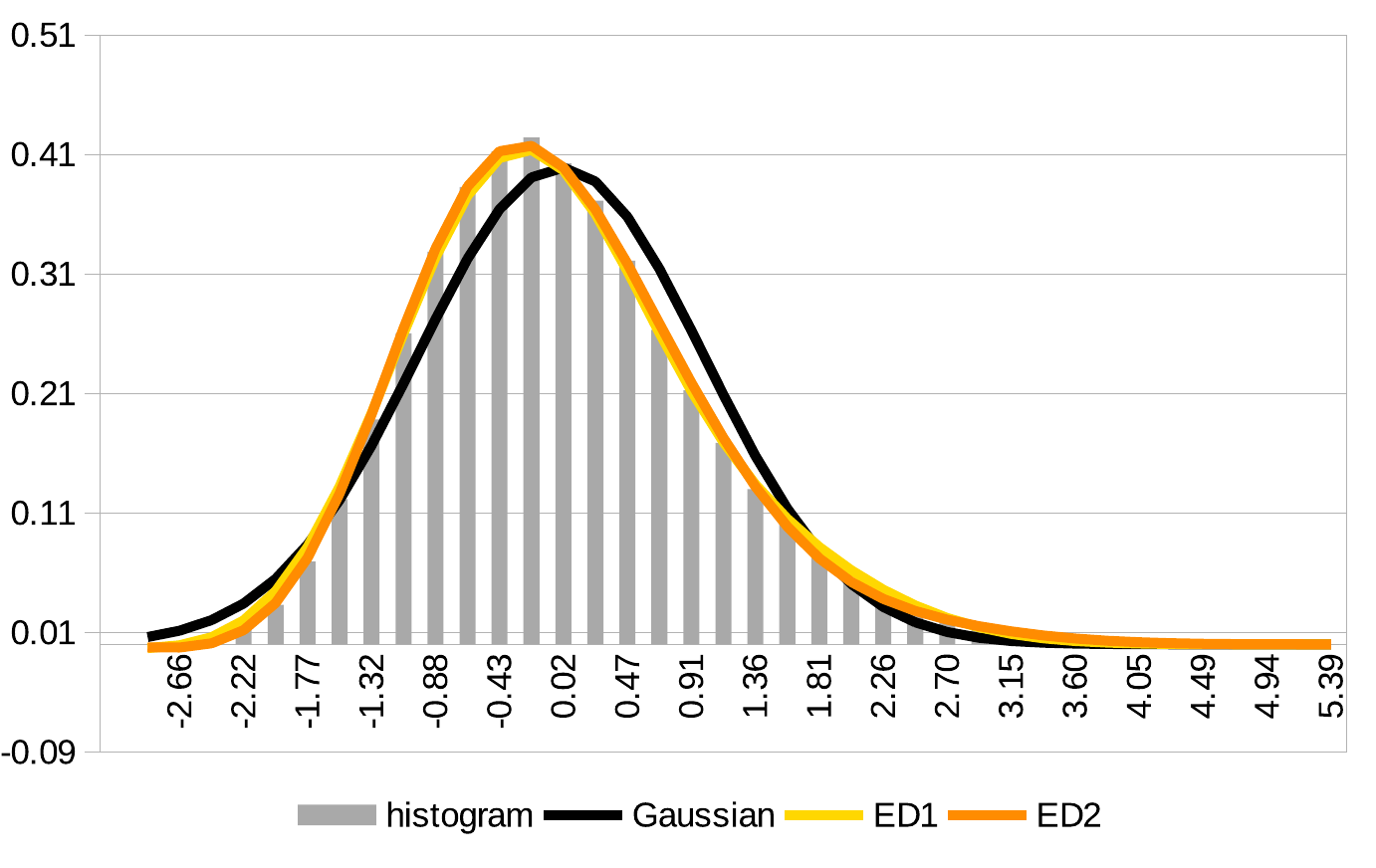} &
    \subfigimg[width=\linewidth]{IV)}{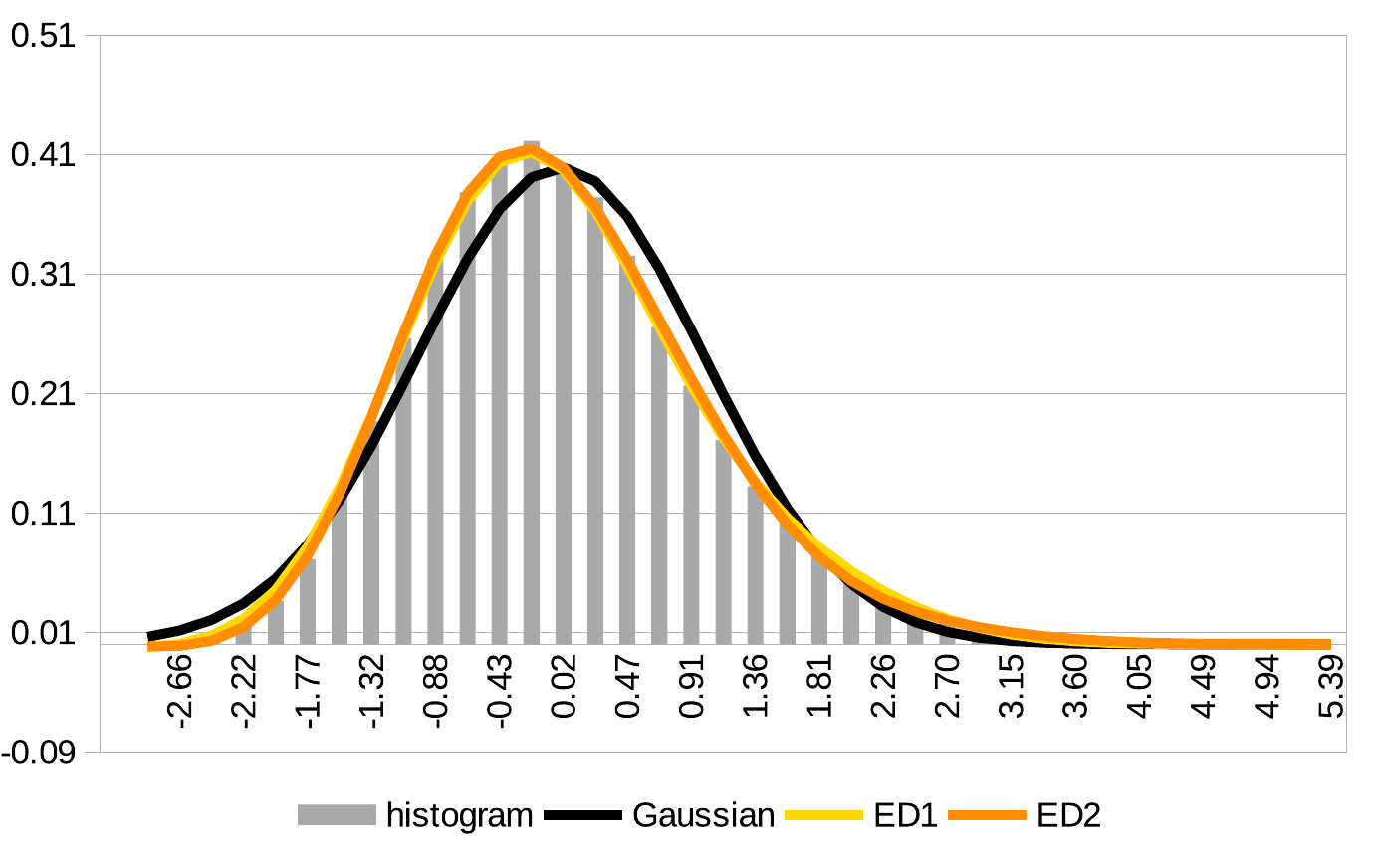} 
    \end{tabular}
\caption{Edgeworth expansion of the standardized QoI in Eq. \eqref{QoI3SGM} for $N=2$, obtained with the SG method.
I): $\sigma_{\gamma}=3$, II): $\sigma_{\gamma}=2.9$, III): $\sigma_{\gamma}=2.8$, IV): $\sigma_{\gamma}=2.7$.}
\label{test2_ED_integral}
\end{figure}

The ED curves are in Figure \ref{test2_ED_integral}. ED2 well approximates the histogram for all values of $\sigma_{\gamma}$, similarly to what was observed in the previous example.
Once again, due to the better approximation provided by ED2 compared to GC3, we conclude that the ED expansion is more valuable than the GC, for the examples considered.

\begin{figure}[!t]
  \centering
 \begin{tabular}{@{}p{0.475\linewidth}@{\quad}p{0.475\linewidth}@{}}
    \subfigimg[width=\linewidth]{I)}{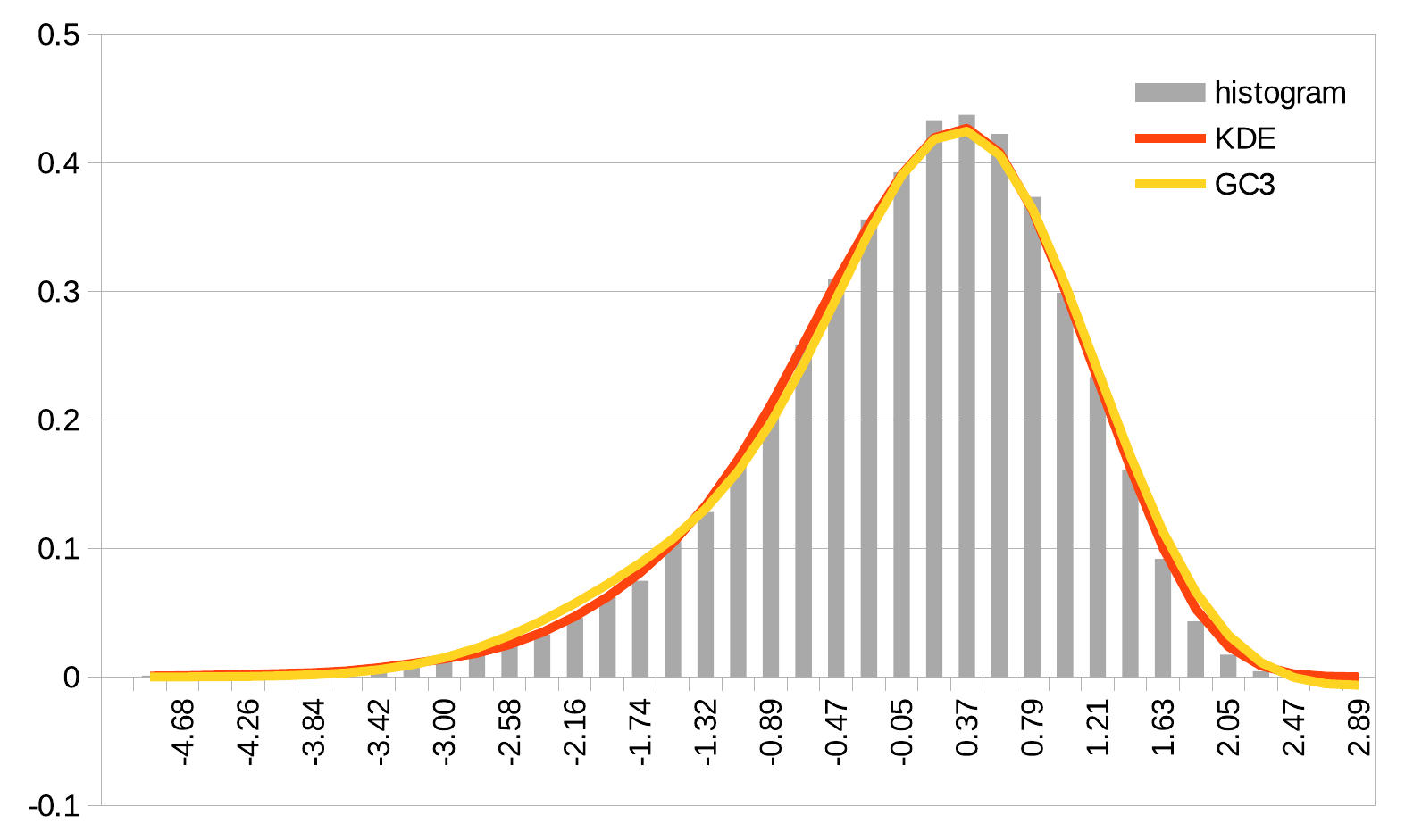}&
    \subfigimg[width=\linewidth]{II)}{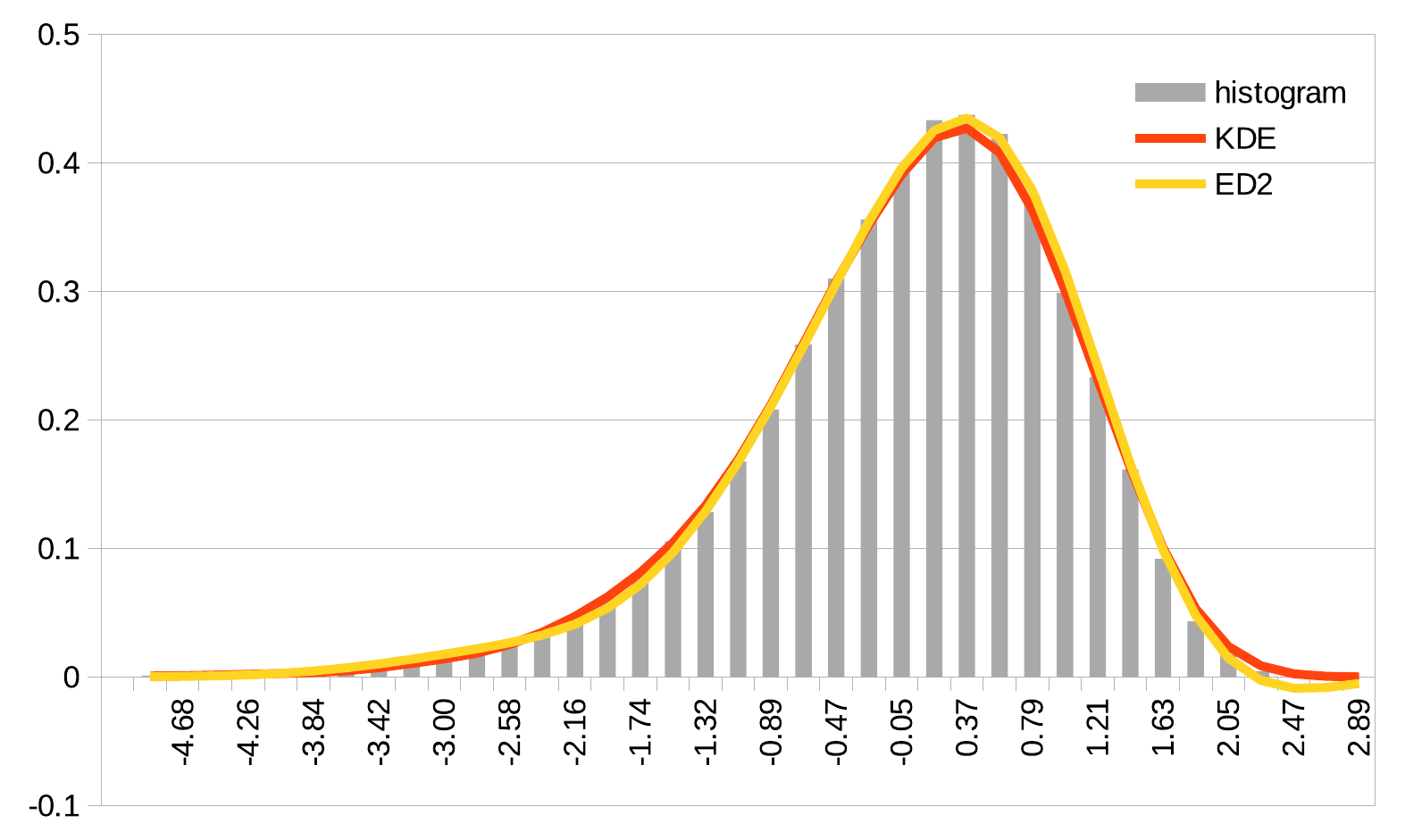} \\
    \subfigimg[width=\linewidth]{III)}{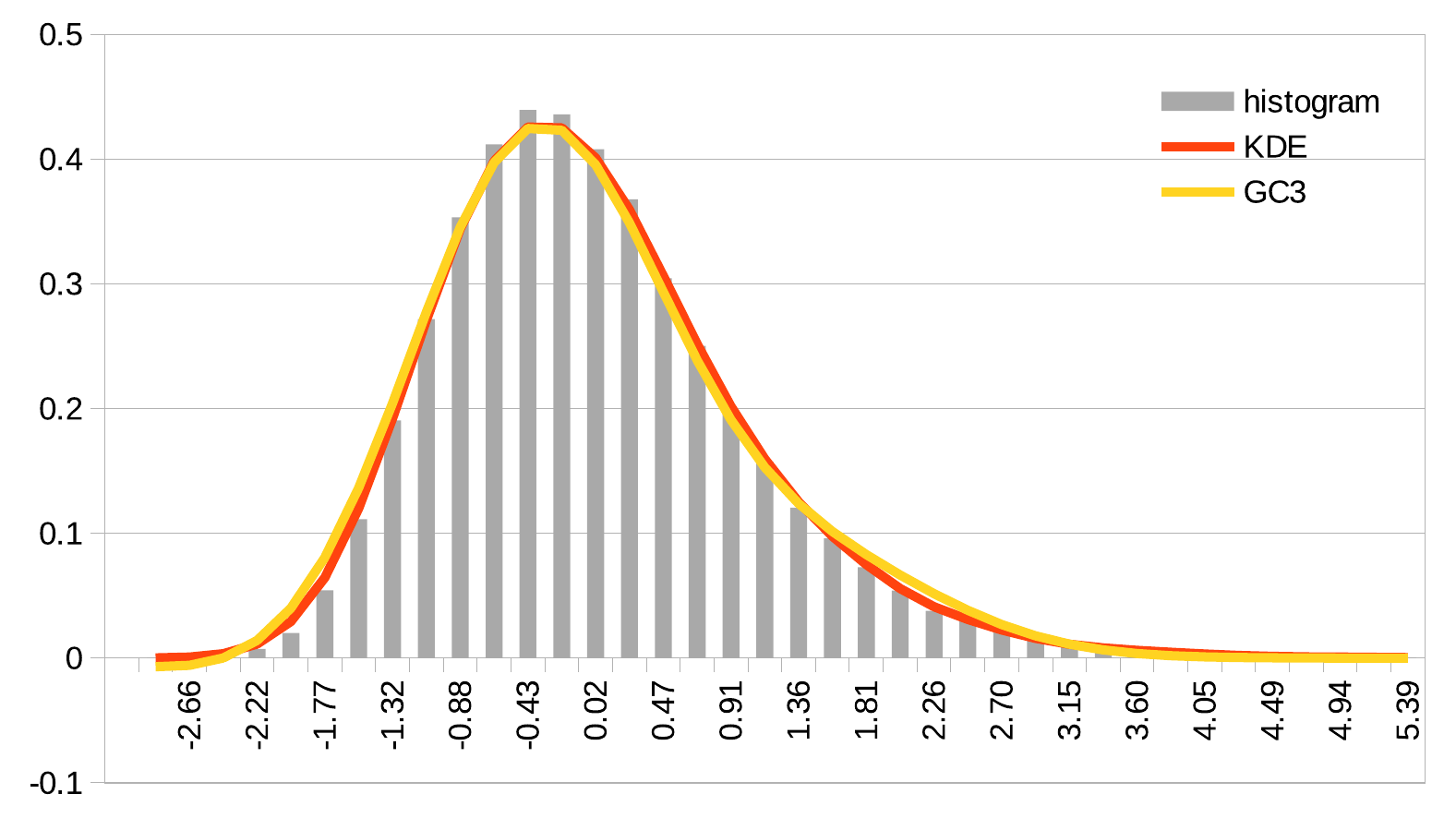} &
    \subfigimg[width=\linewidth]{IV)}{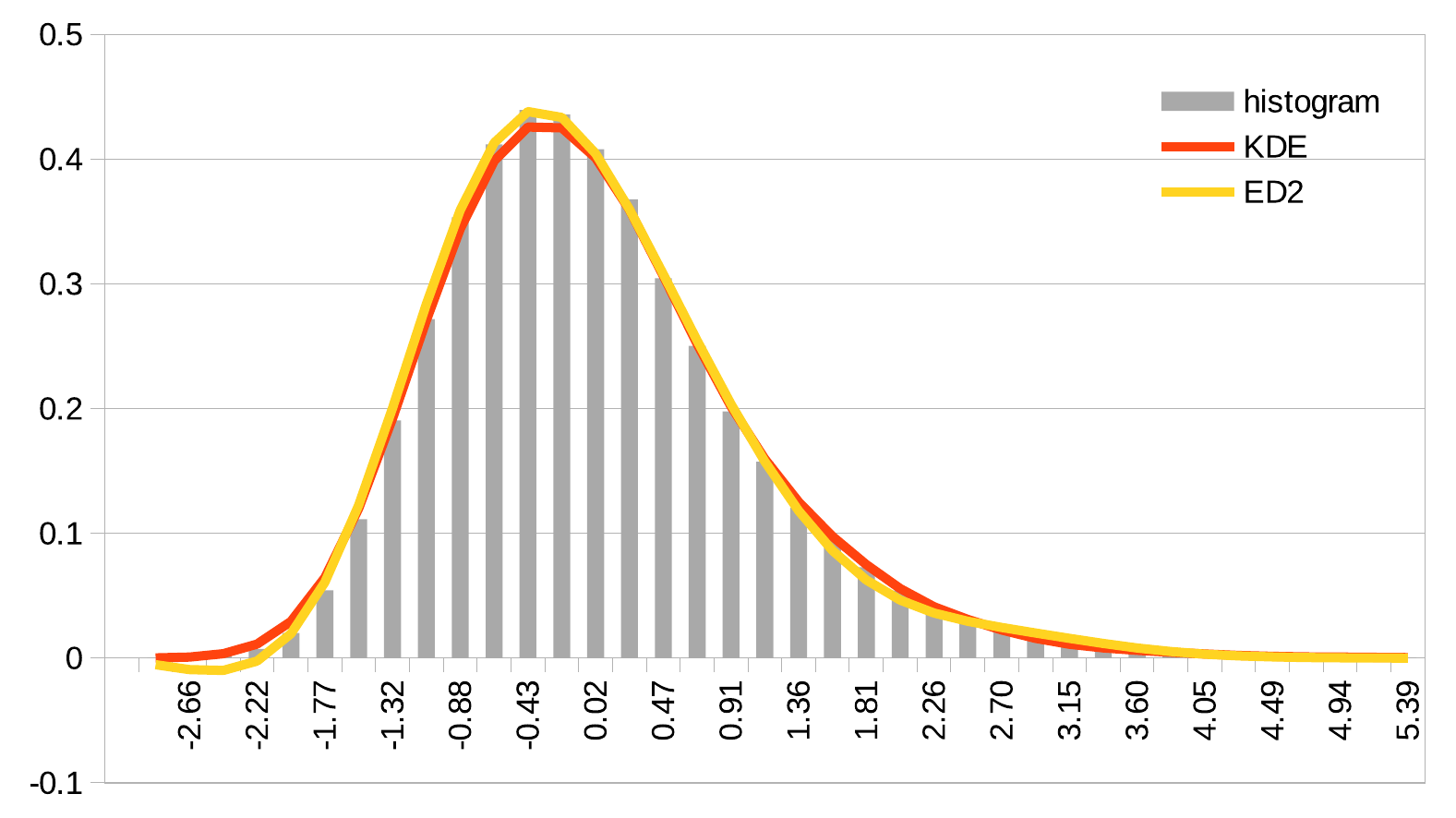} 
    \end{tabular}
\caption{GC3 and ED2 expansions compared with the KDE estimator in Eq. \eqref{KDE}.
I) and II): QoI in \eqref{QoI2SGM} $\sigma_{\gamma}=1.6$, III) and IV): QoI in \eqref{QoI3SGM} $\sigma_{\gamma}=3$.}
\label{KDE_plots}
\end{figure}

\Giacomo{We conclude this section with a comparison of the GC and ED expansions with the kernel density estimator. Given that our analysis is set in a univariate setting, the KDE is given by
\begin{align}\label{KDE}
 f_K(x) = \dfrac{1}{h \, M} \sum\limits_{m=1}^M s\Big(\dfrac{x - \mathcal{Q}_{u}(\bm{\varepsilon}_m)}{h} \Big),
\end{align}
where $s$ is the standard Gaussian distribution defined in Eq.\eqref{univ_pdf}, and $M$ is the size of the sample set.
The parameter $h$ is the bandwidth and we selected it to be the same as the bin width used for the histograms in the previous figures. We chose a standard Gaussian kernel for the KDE because the GC and ED expansions also use a standard Gaussian kernel, hence the comparison is fair.
In Figure \ref{KDE_plots}, GC3 and ED2 are compared to the KDE estimator in Eq. \eqref{KDE} with $M=10^5$:
results for the QoI in \eqref{QoI2SGM} are visible in I) and II) using $\sigma_{\gamma}=1.6$, whereas for the QoI in \eqref{QoI3SGM} they are in III) and IV) considering $\sigma_{\gamma}=3$.
The plots in figure \ref{KDE_plots} show that the GC and ED expansions are comparable to the KDE in terms of accuracy.}
%
%
%
%
%
%
\section{Computational times}
\begin{figure}[h!]
\centering
\includegraphics[scale=0.5]{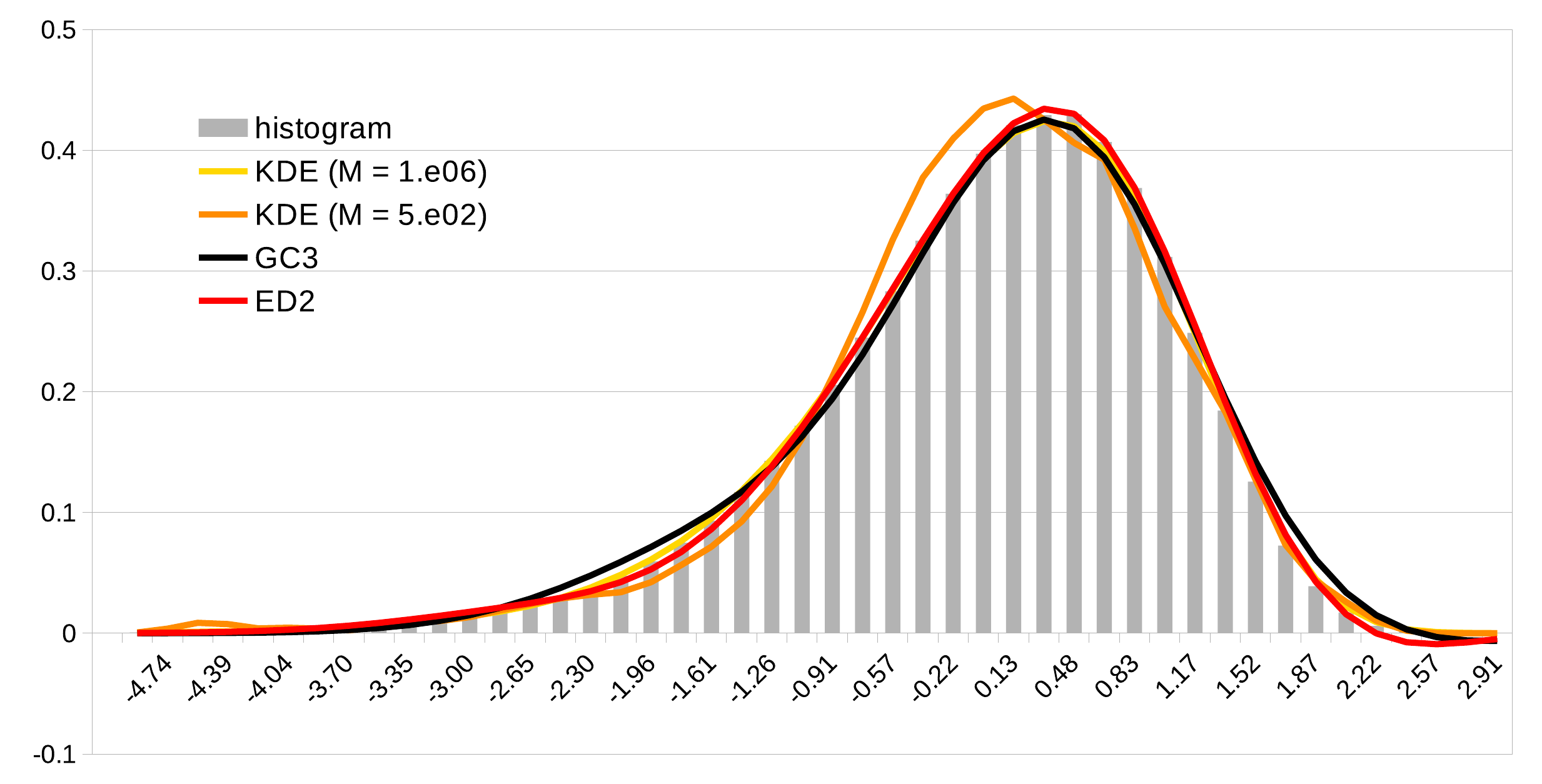}
\caption{KDE and proposed method (GC3 and ED2) curves associated with the results in Table \ref{CPU_times2}.}\label{kde_comp}
\end{figure}

 \Giacomo{Next, we compare the computational time required to construct a histogram approximation using $M=10^6$ samples and the time required by the proposed method. For the latter, the moments in Eq. \eqref{exactMoments} are evaluated with numerical quadrature, and a crude histogram approximation using $M=10^4$ samples is computed.  The comparison is carried out for different values of $N$, i.e. the dimension of the parameter space. The aim of this comparison is to show that the computational time required by the proposed approach is comparable to the one required to obtain an accurate histogram approximation. 
 It is important to remark that the histogram is a discontinuous approximation, regardless of the number of samples employed, whereas the GC and ED expansions are continuous and infinitely differentiable.
 All computations were performed on a Dell Inspiron 15, 5000 series laptop with the CPU \{Intel(R) Core(TM) i3-4030U CPU\@1.90GHz, 1895 MHz\} and 8 GB of RAM.
  The CPU times are shown in Table \ref{CPU_times} and refer to the quantity of interest in Eq. \eqref{QoI2SGM} with $\sigma_{\gamma} = 1.6$. The other simulations parameters are as in the previous section, except $N$ that is varied from 1 to 4. For the proposed method, the total CPU time reported is the sum of three different costs: the cost of the crude approximation with $M=10^4$ samples, the cost of the numerical quadrature that does not involve any sampling, and the cost of performing enough function evaluations to plot the curves, that also does not involve any sampling. 
 \begin {table}[!t]
\setlength\tabcolsep{3.25pt} 
\begin{center}                                                               
	\begin{tabular}{|c|c|c|c|c|c|} \hline	
	   \multicolumn{6}{|c|}{ Computational times [sec] } \\ \cline {1-6}         
	  &  \multicolumn{1}{|c|}{ Histogram } & \multicolumn{4}{|c|}{ Proposed } \\ \cline {2-6}  
 $N$  &  TOT    &  Crude approx & Moments &  Function evals & TOT     \\ \hline 
    1   &  0.843663   & 0.007834  & 0.000115 & 0.000729 & 0.008678   \\ \hline
    2   &  1.251279   & 0.012149  & 0.001191 & 0.000477  & 0.013817    \\ \hline
    3   &  1.788778   & 0.017364  & 0.015390 & 0.000418 &  0.033172    \\ \hline
    4   &  2.533417   & 0.024813  & 0.308295 & 0.000497 &  0.333605    \\ \hline
  	\end{tabular}
\end{center}
\caption{Computational times of histogram and proposed method.}
\label{CPU_times}
\end{table}
The results show that, for the values of $N$ considered, the CPU times of the proposed method are lower than those of the histogram. It is true however that as the values of $N$ increase, the number of quadrature points necessary for the evaluation of the moments in Eq. \eqref{exactMoments} will grow, likely causing the histogram to eventually become faster.
In Table \ref{CPU_times2} we compare the CPU time of a function evaluation using the KDE with the CPU time of a function evaluation using the proposed method.
The main difference between the two methods is that the dominant cost for the KDE is associated with online operations, whereas for the proposed method the most expensive operations are done offline.
For the KDE, every function evaluation requires the computation of as many kernel values as the total number of samples.
Function evaluations count as online operations for the estimation of the PDF.
If the sample set is composed of $M_s$ points and $M_e$ function evaluations are required to estimate the PDF, the estimated cost of the KDE is proportional to the product $M_s \, M_e$ (the influence of the parameter space dimension is negligible).
For the proposed method, the cost of a single function evaluation is given by the computation of the moments with numerical quadrature and by the crude histogram computation. These operations are all offline and can be performed once and for all regardless of the number of functions evaluations necessary to estimate the PDF.
If $M_c < M_s$ is the number of samples used for the crude histogram and $Q_p$ is the number of quadrature points required for the exact one-dimensional numerical quadrature, the estimated cost of the proposed method is proportional to $M_c + Q_p^N$, with $N$ being the dimension of the parameter space.
Table \ref{CPU_times2} shows CPU times results for $N=1,2,3,4$, considering $M_c = 10^4$ and $M_s = 10^6$ or $M_s = 500$ for the KDE. The KDE with $M_s = 500$ is the fastest approximation, however it also the most inaccurate, as it can be seen from Figure \ref{kde_comp}. The KDE with $M_s = 10^6$, GC3 and ED2 have a comparable level of accuracy, hence we focus our analysis on the comparison between the proposed method and the KDE with $M_s = 10^6$. For $N=1,2,3$ the proposed method is faster than a single function evaluation of the KDE with $M_s = 10^6$, whereas for $N=4$, one KDE evaluation is approximately nine times faster than the proposed method. It is very likely that in general more than just nine function evaluations will be necessary to appropriately describe the approximated PDF. For instance, in Figure \ref{kde_comp}, the plots have been obtained with $46$ function evaluations. Hence, the proposed method can compete with the KDE in terms of CPU time.}
 \begin {table}[!t]
\setlength\tabcolsep{3.25pt} 
\begin{center}                                                               
	\begin{tabular}{|c|c|c|c|c|c|c|} \hline	
	   \multicolumn{7}{|c|}{ CPU times for a single function evaluation [sec] } \\ \cline {1-7} 
	 &   \multicolumn{2}{|c|}{ KDE ($M=10^6$) } & \multicolumn{2}{|c|}{ KDE ($M=500$)  } & \multicolumn{2}{|c|}{ Proposed } \\ \cline {2-7}  
	$N$ &  Offline & Online &  Offline & Online & Offline & Online \\ \hline 
	1    & negligible & 0.048431  & negligible  & 0.000034  & 0.007968    &  negligible  \\ \hline 
        2    & negligible & 0.047788  & negligible  & 0.000036  & 0.012642   &  negligible  \\ \hline 
        3    & negligible & 0.052763  & negligible  & 0.000050  & 0.044448   &  negligible  \\ \hline 
        4    & negligible & 0.047178  & negligible  & 0.000037  & 0.334770   &  negligible  \\ \hline 
  	\end{tabular}
\end{center}
\caption{Computational times for a single evaluation of the KDE and of the proposed method ($10^4$ are used for the crude histogram).}
\label{CPU_times2}
\end{table}
\section{Discussion}\label{discussion}
It has been shown that the GC and ED truncated expansions represent a valid alternative to existing methods for the approximation of probability density functions associated with solutions of {PDEs with random parameters}.
Our numerical results suggested that GC and ED provide an accurate estimate when the PDF is nearly Gaussian.
This is consistent with the nature of the distributions.
Even in the case of a {non-Gaussian} PDF, the truncated expansions well approximated the distributions. The asymptotic character of the ED expansion makes it more valuable than the GC due to a better monitoring of the error, as the truncation order is increased.
Moreover, ED better approximated the distribution than GC, at least up to the truncation order considered in this work.

The proposed method is easy to implement and takes advantage of the fact that exact moments can be computed with the {SG method}, provided that enough quadrature points are considered.
\Giacomo{Moreover, all the computational burden is associated with offline computations, and point-wise evaluations have a negligible cost}.
A limitation is the lack of a rigorous procedure to determine the necessary number of terms in the truncated series to obtain the best possible approximation. This is inherent in the use of truncated series and the issue is present in all works on ED and GC that we were able to find in the literature.  While the optimal number of terms to retain will likely be too large for the desired error in case of a convergent series, it could be theoretically determined a priori in case of an asymptotic expansion such as Edgeworth. Unfortunately, the lack of explicit error bounds make this unfeasible at the moment.
\section*{Conflict of interest}
The authors declare that they have no conflict of interest.

\bibliographystyle{spmpsci}
\bibliography{expansion}

\begin{thebibliography}{10}
\providecommand{\url}[1]{{#1}}
\providecommand{\urlprefix}{URL }
\expandafter\ifx\csname urlstyle\endcsname\relax
  \providecommand{\doi}[1]{DOI~\discretionary{}{}{}#1}\else
  \providecommand{\doi}{DOI~\discretionary{}{}{}\begingroup
  \urlstyle{rm}\Url}\fi

\bibitem{aulisa2018construction}
Aulisa, E., Capodaglio, G., Ke, G.: Construction of h-refined continuous finite
  element spaces with arbitrary hanging node configurations and applications to
  multigrid algorithms.
\newblock arXiv preprint arXiv:1804.10632  (2018)

\bibitem{babuvska2007stochastic}
Babu{\v{s}}ka, I., Nobile, F., Tempone, R.: A stochastic collocation method for
  elliptic partial differential equations with random input data.
\newblock SIAM Journal on Numerical Analysis \textbf{45}(3), 1005--1034 (2007)

\bibitem{babuska2004galerkin}
Babuska, I., Tempone, R., Zouraris, G.E.: Galerkin finite element
  approximations of stochastic elliptic partial differential equations.
\newblock SIAM Journal on Numerical Analysis \textbf{42}(2), 800--825 (2004)

\bibitem{babuvska2005solving}
Babu{\v{s}}ka, I., Tempone, R., Zouraris, G.E.: Solving elliptic boundary value
  problems with uncertain coefficients by the finite element method: the
  stochastic formulation.
\newblock Computer methods in applied mechanics and engineering
  \textbf{194}(12-16), 1251--1294 (2005)

\bibitem{bell1927partition}
Bell, E.T.: Partition polynomials.
\newblock Annals of Mathematics pp. 38--46 (1927)

\bibitem{berberan2007expressing}
Berberan-Santos, M.N.: Expressing a probability density function in terms of
  another {PDF}: A generalized {G}ram-{C}harlier expansion.
\newblock Journal of Mathematical Chemistry \textbf{42}(3), 585--594 (2007)

\bibitem{blinnikov1998expansions}
Blinnikov, S., Moessner, R.: Expansions for nearly {G}aussian distributions.
\newblock Astronomy and Astrophysics Supplement Series \textbf{130}(1),
  193--205 (1998)

\bibitem{brenn2017revisit}
Brenn, T., Anfinsen, S.N.: A revisit of the {G}ram-{C}harlier and {E}dgeworth
  series expansions  (2017)

\bibitem{brenner2007mathematical}
Brenner, S., Scott, R.: The mathematical theory of finite element methods,
  vol.~15.
\newblock Springer Science \& Business Media (2007)

\bibitem{bui2013computational}
Bui-Thanh, T., Ghattas, O., Martin, J., Stadler, G.: A computational framework
  for infinite-dimensional {B}ayesian inverse problems part i: The linearized
  case, with application to global seismic inversion.
\newblock SIAM Journal on Scientific Computing \textbf{35}(6), A2494--A2523
  (2013)

\bibitem{femus-web-page}
Capodaglio, G.: Github webpage.
\newblock \urlprefix\url{https://github.com/gcapodag/MyFEMuS}

\bibitem{chacon2013data}
Chac{\'o}n, J.E., Duong, T., et~al.: Data-driven density derivative estimation,
  with applications to nonparametric clustering and bump hunting.
\newblock Electronic Journal of Statistics \textbf{7}, 499--532 (2013)

\bibitem{chen2016model}
Chen, P., Schwab, C.: Model order reduction methods in computational
  uncertainty quantification.
\newblock Handbook of Uncertainty Quantification pp. 1--53 (2016)

\bibitem{cheng1995mean}
Cheng, Y.: Mean shift, mode seeking, and clustering.
\newblock IEEE transactions on pattern analysis and machine intelligence
  \textbf{17}(8), 790--799 (1995)

\bibitem{ciarlet2002finite}
Ciarlet, P.G.: The finite element method for elliptic problems, vol.~40.
\newblock Siam (2002)

\bibitem{cliffe2011multilevel}
Cliffe, K.A., Giles, M.B., Scheichl, R., Teckentrup, A.L.: Multilevel {M}onte
  {C}arlo methods and applications to elliptic {PDE}s with random coefficients.
\newblock Computing and Visualization in Science \textbf{14}(1), 3 (2011)

\bibitem{contaldi1999photographing}
Contaldi, C.R., Bean, R., Magueijo, J.: Photographing the wave function of the
  universe.
\newblock Physics Letters B \textbf{468}(3-4), 189--194 (1999)

\bibitem{cramer1928composition}
Cram{\'e}r, H.: On the composition of elementary errors: First paper:
  Mathematical deductions.
\newblock Scandinavian Actuarial Journal \textbf{1928}(1), 13--74 (1928)

\bibitem{cramer2016mathematical}
Cram{\'e}r, H.: Mathematical methods of statistics (PMS-9), vol.~9.
\newblock Princeton university press (2016)

\bibitem{dashti2017bayesian}
Dashti, M., Stuart, A.M.: The {B}ayesian approach to inverse problems.
\newblock Handbook of Uncertainty Quantification pp. 311--428 (2017)

\bibitem{de2011edgeworth}
De~Kock, M., Eggers, H., Schmiegel, J.: {E}dgeworth versus {G}ram-{C}harlier
  series: x-cumulant and probability density tests.
\newblock Physics of Particles and Nuclei Letters \textbf{8}(9), 1023--1027
  (2011)

\bibitem{di2001mathematical}
Di~Marco, V.B., Bombi, G.G.: Mathematical functions for the representation of
  chromatographic peaks.
\newblock Journal of Chromatography A \textbf{931}(1-2), 1--30 (2001)

\bibitem{eggers2011determining}
Eggers, H.C., de~Kock, M.B., Schmiegel, J.: Determining source cumulants in
  femtoscopy with {G}ram--{C}harlier and {E}dgeworth series.
\newblock Modern Physics Letters A \textbf{26}(24), 1771--1782 (2011)

\bibitem{fan2012probabilistic}
Fan, M., Vittal, V., Heydt, G.T., Ayyanar, R.: Probabilistic power flow studies
  for transmission systems with photovoltaic generation using cumulants.
\newblock IEEE Transactions on Power Systems \textbf{27}(4), 2251--2261 (2012)

\bibitem{fishman2013monte}
Fishman, G.: {M}onte {C}arlo: concepts, algorithms, and applications.
\newblock Springer Science \& Business Media (2013)

\bibitem{frauenfelder2005finite}
Frauenfelder, P., Schwab, C., Todor, R.A.: Finite elements for elliptic
  problems with stochastic coefficients.
\newblock Computer methods in applied mechanics and engineering
  \textbf{194}(2-5), 205--228 (2005)

\bibitem{fukunaga1975estimation}
Fukunaga, K., Hostetler, L.: The estimation of the gradient of a density
  function, with applications in pattern recognition.
\newblock IEEE Transactions on information theory \textbf{21}(1), 32--40 (1975)

\bibitem{ghanem1991stochastic}
Ghanem, R.G., Spanos, P.D.: Stochastic finite element method: Response
  statistics.
\newblock In: Stochastic Finite Elements: A Spectral Approach, pp. 101--119.
  Springer (1991)

\bibitem{gunzburger2014stochastic}
Gunzburger, M.D., Webster, C.G., Zhang, G.: Stochastic finite element methods
  for partial differential equations with random input data.
\newblock Acta Numerica \textbf{23}, 521--650 (2014)

\bibitem{hernandez2005slepc}
Hernandez, V., Roman, J.E., Vidal, V.: {SLEP}c: A scalable and flexible toolkit
  for the solution of eigenvalue problems.
\newblock ACM Transactions on Mathematical Software (TOMS) \textbf{31}(3),
  351--362 (2005)

\bibitem{huang2001convergence}
Huang, S., Quek, S., Phoon, K.: Convergence study of the truncated
  {K}arhunen--{L}oeve expansion for simulation of stochastic processes.
\newblock International journal for numerical methods in engineering
  \textbf{52}(9), 1029--1043 (2001)

\bibitem{izenman1991review}
Izenman, A.J.: Review papers: Recent developments in nonparametric density
  estimation.
\newblock Journal of the American Statistical Association \textbf{86}(413),
  205--224 (1991)

\bibitem{jondeau2001gram}
Jondeau, E., Rockinger, M.: {G}ram--{C}harlier densities.
\newblock Journal of Economic Dynamics and Control \textbf{25}(10), 1457--1483
  (2001)

\bibitem{juszkiewicz1993weakly}
Juszkiewicz, R., Weinberg, D., Amsterdamski, P., Chodorowski, M., Bouchet, F.:
  Weakly non-linear {G}aussian fluctuations and the {E}dgeworth expansion.
\newblock arXiv preprint astro-ph/9308012  (1993)

\bibitem{kendall1943advanced}
Kendall, M.G.: Advanced Theory Of Statistics Vol-I.
\newblock Charles Griffin: London (1943)

\bibitem{li2008fourier}
Li, C., Feng, Y., Owen, D., Li, D., Davis, I.: A
  {F}ourier--{K}arhunen--{L}o{\`e}ve discretization scheme for stationary
  random material properties in {SFEM}.
\newblock International journal for numerical methods in engineering
  \textbf{73}(13), 1942--1965 (2008)

\bibitem{ma2009efficient}
Ma, X., Zabaras, N.: An efficient {B}ayesian inference approach to inverse
  problems based on an adaptive sparse grid collocation method.
\newblock Inverse Problems \textbf{25}(3), 035013 (2009)

\bibitem{marzouk2007stochastic}
Marzouk, Y.M., Najm, H.N., Rahn, L.A.: Stochastic spectral methods for
  efficient {B}ayesian solution of inverse problems.
\newblock Journal of Computational Physics \textbf{224}(2), 560--586 (2007)

\bibitem{metropolis1949monte}
Metropolis, N., Ulam, S.: The {M}onte {C}arlo method.
\newblock Journal of the American statistical association \textbf{44}(247),
  335--341 (1949)

\bibitem{mihoubi2008bell}
Mihoubi, M.: Bell polynomials and binomial type sequences.
\newblock Discrete Mathematics \textbf{308}(12), 2450--2459 (2008)

\bibitem{niguez2012forecasting}
{\~N}{\'\i}guez, T.M., Perote, J.: Forecasting heavy-tailed densities with
  positive {E}dgeworth and {G}ram-{C}harlier expansions.
\newblock Oxford Bulletin of Economics and Statistics \textbf{74}(4), 600--627
  (2012)

\bibitem{nobile2008anisotropic}
Nobile, F., Tempone, R., Webster, C.G.: An anisotropic sparse grid stochastic
  collocation method for partial differential equations with random input data.
\newblock SIAM Journal on Numerical Analysis \textbf{46}(5), 2411--2442 (2008)

\bibitem{nobile2008sparse}
Nobile, F., Tempone, R., Webster, C.G.: A sparse grid stochastic collocation
  method for partial differential equations with random input data.
\newblock SIAM Journal on Numerical Analysis \textbf{46}(5), 2309--2345 (2008)

\bibitem{noh2014bias}
Noh, Y.K., Sugiyama, M., Liu, S., Plessis, M.C., Park, F.C., Lee, D.D.: Bias
  reduction and metric learning for nearest-neighbor estimation of
  {K}ullback-{L}eibler divergence.
\newblock In: Artificial Intelligence and Statistics, pp. 669--677 (2014)

\bibitem{o1992using}
O'brien, M.: Using the {G}ram-{C}harlier expansion to produce vibronic band
  shapes in strong coupling.
\newblock Journal of Physics: Condensed Matter \textbf{4}(9), 2347 (1992)

\bibitem{olive1991gram}
Oliv{\'e}, J., Grimalt, J.O.: {G}ram-{C}harlier and {E}dgeworth-{C}ram{\'e}r
  series in the characterization of chromatographic peaks.
\newblock Analytica chimica acta \textbf{249}(2), 337--348 (1991)

\bibitem{pender2014Gram}
Pender, J.: {G}ram {C}harlier expansion for time varying multiserver queues
  with abandonment.
\newblock SIAM Journal on Applied Mathematics \textbf{74}(4), 1238--1265 (2014)

\bibitem{petrov1995limit}
Petrov, V.V.: Limit theorems of probability theory: sequences of independent
  random variables.
\newblock Tech. rep., Oxford, New York (1995)

\bibitem{petrov2012sums}
Petrov, V.V.: Sums of independent random variables, vol.~82.
\newblock Springer Science \& Business Media (2012)

\bibitem{popovic2012easy}
Popovic, R., Goldsman, D.: Easy {G}ram-{C}harlier valuations of options.
\newblock Journal of Derivatives \textbf{20}(2), 79 (2012)

\bibitem{rickman2015calculating}
Rickman, J., Lawrence, A., Rollett, A., Harmer, M.: Calculating probability
  densities associated with grain-size distributions.
\newblock Computational Materials Science \textbf{101}, 211--215 (2015)

\bibitem{sasaki2015direct}
Sasaki, H., Noh, Y.K., Sugiyama, M.: Direct density-derivative estimation and
  its application in {KL}-divergence approximation.
\newblock In: Artificial Intelligence and Statistics, pp. 809--818 (2015)

\bibitem{schevenels2004application}
Schevenels, M., Lombaert, G., Degrande, G.: Application of the stochastic
  finite element method for gaussian and non-gaussian systems.
\newblock In: ISMA2004 International Conference on Noise and Vibration
  Engineering, pp. 3299--3314 (2004)

\bibitem{sedgewick2013introduction}
Sedgewick, R., Flajolet, P.: An introduction to the analysis of algorithms.
\newblock Pearson Education India (2013)

\bibitem{stuart2010inverse}
Stuart, A.M.: Inverse problems: a {B}ayesian perspective.
\newblock Acta Numerica \textbf{19}, 451--559 (2010)

\bibitem{tartakovsky2011pdf}
Tartakovsky, D.M., Broyda, S.: Pdf equations for advective--reactive transport
  in heterogeneous porous media with uncertain properties.
\newblock Journal of contaminant hydrology \textbf{120}, 129--140 (2011)

\bibitem{wallace1958asymptotic}
Wallace, D.L.: Asymptotic approximations to distributions.
\newblock The Annals of Mathematical Statistics \textbf{29}(3), 635--654 (1958)

\bibitem{wan2005adaptive}
Wan, X., Karniadakis, G.E.: An adaptive multi-element generalized polynomial
  chaos method for stochastic differential equations.
\newblock Journal of Computational Physics \textbf{209}(2), 617--642 (2005)

\bibitem{xiu2005high}
Xiu, D., Hesthaven, J.S.: High-order collocation methods for differential
  equations with random inputs.
\newblock SIAM Journal on Scientific Computing \textbf{27}(3), 1118--1139
  (2005)

\bibitem{zapevalov2011simulating}
Zapevalov, A., Bol’shakov, A., Smolov, V.: Simulating of the probability
  density of sea surface elevations using the {G}ram-{C}harlier series.
\newblock Oceanology \textbf{51}(3), 407--414 (2011)

\end{thebibliography}

\end{document}